\documentclass[11pt,a4paper,twoside,enabledeprecatedfontcommands,final]{scrartcl}
\usepackage{a4wide}
\usepackage{amsfonts}
\usepackage{amsmath}
\usepackage{amssymb}
\usepackage{amsthm}
\usepackage[numbers]{natbib}
\usepackage{booktabs}
\usepackage{algorithm,algpseudocode}
\usepackage{pgfplots}
\pgfplotsset{compat=1.8}
\usepackage{graphicx}
\usepackage{color}
\usepackage{colortbl}   
\usepackage{showkeys}
\usepackage{scrpage2}
\usepackage{multirow}
\usepackage{boxedminipage}
\usepackage{enumerate}
\usepackage{subfig}
\usepackage{textcomp}

\usepackage{hyperref}

\usepackage{todonotes}
\usepackage{placeins}

\usepackage{siunitx}
\newcommand{\N}{\ensuremath{\mathbb{N}}}

\newcommand{\T}{\ensuremath{\mathbb{T}}}

\newcommand{\Z}{\ensuremath{\mathbb{Z}}}

\newcommand{\R}{\ensuremath{\mathbb{R}}}
\newcommand{\C}{\ensuremath{\mathbb{C}}}

\newcommand{\ii}{\textnormal{i}}
\newcommand{\e}{\textnormal{e}}

\newcommand{\ceil}[1]{\left\lceil#1\right\rceil}
\newcommand{\floor}[1]{\left\lfloor#1\right\rfloor}
\newcommand{\zb}[1]{\ensuremath{\boldsymbol{#1}}}

\renewcommand{\ln}{\mathrm{ln\,}}
\newcommand{\bigtimes}{\mathop{\text{\Large{$\times$}}}}

\newcommand{\boldx}{{\ensuremath{\boldsymbol{x}}}}
\newcommand{\boldA}{{\ensuremath{\boldsymbol{A}}}}
\newcommand{\boldy}{{\ensuremath{\boldsymbol{y}}}}
\newcommand{\boldk}{{\ensuremath{\boldsymbol{k}}}}
\newcommand{\boldh}{{\ensuremath{\boldsymbol{h}}}}

\newcommand{\boldz}{{\ensuremath{\boldsymbol{z}}}}

\newcommand{\boldj}{{\ensuremath{\boldsymbol{j}}}}

\newcommand{\boldone}{{\ensuremath{\boldsymbol{1}}}}

\definecolor{darkgreen}{rgb}{0.0,0.5,0.0}

\newtheorem{theorem}{Theorem}[section]
\newtheorem{lemma}[theorem]{Lemma}
\newtheorem{remark}[theorem]{Remark}
\newtheorem{generalisation}[theorem]{Generalisation}
\newtheorem{definition}[theorem]{Definition}
\newtheorem{example}[theorem]{Example}
\newtheorem{corollary}[theorem]{Corollary}
\newtheorem{proposition}[theorem]{Proposition}

\newenvironment{Theorem}{\goodbreak \begin{theorem}\sl}{\end{theorem}}
\newenvironment{Lemma}{\goodbreak \begin{lemma}\sl}{\end{lemma}}
\newenvironment{Remark}{\goodbreak \begin{remark}\rm}{\bend\end{remark}}

\newenvironment{Corollary}{\goodbreak \begin{corollary}\sl}{\end{corollary}}

\algnewcommand\algorithmicforeach{\textbf{for each}}
\algdef{S}[FOR]{ForEach}[1]{\algorithmicforeach\ #1\ \algorithmicdo}

\makeatletter
\def\imod#1{\allowbreak\mkern10mu({\operator@font mod}\,\,#1)}
\makeatother

\numberwithin{equation}{section}
\numberwithin{table}{section}
\numberwithin{figure}{section}

\newcommand{\bend}{\hspace*{0ex} \hfill \hbox{\vrule height
    1.5ex\vbox{\hrule width 1.4ex \vskip 1.4ex\hrule  width 1.4ex}\vrule
    height 1.5ex}}

\long\def\symbolfootnote[#1]#2{\begingroup%
\def\thefootnote{\fnsymbol{footnote}}\footnote[#1]{#2}\endgroup}

\newcommand{\OO}[1]{\mathcal{O}\left(#1\right)}

\renewcommand{\mathbf}[1]{\ensuremath{\boldsymbol{#1}}}
\renewcommand{\textbf}[1]{{\ensuremath{\boldsymbol{#1}}}}
\renewcommand{\thefootnote}{\fnsymbol{footnote}}

\title{A sample efficient sparse FFT for arbitrary frequency candidate sets in high dimensions}
\date{}
\author{Lutz K\"ammerer\footnotemark[1] \and Felix Krahmer\footnotemark[2] \and 
        Toni Volkmer\footnotemark[3]}
        
        \hypersetup{
          pdfcreator = {pdflatex},
          plainpages=false,
          pdfstartview=Fit,       
          pdfpagemode=UseOutlines, 
          bookmarksnumbered=true, 
          bookmarksopen=false,     
          bookmarksopenlevel=0,   
          colorlinks=true,       
          linkcolor=black,
          citecolor=black,
          urlcolor=black}

\begin{document}

\maketitle

\footnotetext[1]{
  Chemnitz University of Technology, Faculty of Mathematics, 09107 Chemnitz, Germany\\
  lutz.kaemmerer@math.tu-chemnitz.de, Phone:+49-371-531-37728, Fax:+49-371-531-837728
}
\footnotetext[2]{
  Department of Mathematics, Technische Universit\"at M\"unchen, 85748 Garching/Munich, Germany \\
  felix.krahmer@tum.de
}
\footnotetext[3]{
  Chemnitz University of Technology, Faculty of Mathematics, 09107 Chemnitz, Germany\\
  toni.volkmer@math.tu-chemnitz.de, Phone:+49-371-531-39999, Fax:+49-371-531-839999
}

\begin{abstract}

In this paper a sublinear time algorithm is presented for the reconstruction of functions that can be represented by just few out of a potentially large candidate set of Fourier basis functions in high spatial dimensions, a so-called high-dimensional sparse fast Fourier transform. In contrast to many other such algorithms, our method works for arbitrary candidate sets and does not make additional structural assumptions on the candidate set. Our transform  significantly improves upon the other approaches available for such a general framework in terms of the scaling of the sample complexity. Our algorithm is based on sampling the function along multiple rank-1 lattices with random generators.

Combined with a dimension-incremental approach, our method yields a sparse Fourier transform whose computational complexity only grows mildly in the dimension and can hence be efficiently computed even in high dimensions.
Our theoretical analysis establishes that any Fourier $s$-sparse function can be accurately reconstructed with high probability. This guarantee is complemented by several numerical tests demonstrating 
the high efficiency and versatile applicability for the exactly sparse case and also for the compressible case.

\medskip

\medskip

\end{abstract}
\newpage

\section{Introduction}
The problem of reconstructing or approximating a function from a finite number of samples is one of the central objects of study in approximation theory. In finite dimensional function spaces (or function spaces that allow for good finite-dimensional approximations), this problem can be formulated as a linear system $\boldA \boldx = \boldy$. Here $\boldy$ is the vector of observed function values, $\boldx$ is the (unknown) coefficient vector of the function in a given basis $B$, and the matrix~$\boldA$ represents the linear map that maps basis coefficients to function evaluations; knowing the basis $B$ and the sampling points, $\boldA$ can be explicitly determined. The focus of this paper will be on the case where $\boldA$ is the discrete  Fourier transform (DFT) matrix or some submatrix of it.

When $\boldA$ is left-invertible with left inverse $\boldsymbol{L}$, as it is the case for the full DFT matrix, one can find $\boldx$ via the the matrix-vector product of $\boldsymbol{L}$ and $\boldy$. The computational complexity of such a matrix-vector multiplication can, in general, not be improved beyond a linear scaling in the number of entries of the matrix even if the matrix inverse has been precomputed, as this is what is needed just to read the matrix.
In the case of DFT matrices, further acceleration is possible via the  
Fast Fourier transform (FFT)~\cite{CoTu65}, which improves the complexity of multiplying the (inverse) discrete Fourier matrix of size  $n\times n$ with an arbitrary vector of length $n$
from $\mathcal{O}(n^2)$ to $\mathcal{O}(n\log n)$. An extension of this result for irregular sampling nodes with similar advantages is the nonequispaced fast Fourier transform (see e.g.\ \cite{PoStTa01}). Further optimization is possible if the frequencies, i.e., the multi-indices of the basis functions, belong to special structured sets such as hyperbolic crosses, cf.~\cite{DoKuPo10}.

A natural next step is to consider functions where the frequencies of the non-zero Fourier coefficients form a relatively small set $I\subset\Z^d$ within a possibly huge candidate set $\Gamma\subset\Z^d$, but for which the set $I$ is unknown. Such functions are commonly referred to as Fourier sparse. Here a function is $s$-sparse in the given basis if it can be expressed as a linear combination of only at most $s$ basis elements, i.e., $s\ge |I|$. Naturally, this reduced model approach is assumed in order to reduce both the number of samples required and the computation time. The crucial challenge in this approach is the identification of $I$, which is a much harder task compared to the classical approach where only the Fourier coefficients, i.e., the entries of coefficient vector~$\boldx$, are unknown.

The focus in the area of compressive sensing  has been on the sample complexity of this problem. Here, one seeks to recover a sparse signal from as few samples as possible.
It has been shown that recovery is possible for a number of random samples that scales linearly in the sparsity~$s$ up to logarithmic factors \cite{CaRoTa06, RuVe08}, even when the signal is sparse in an arbitrary incoherent basis~\cite{CaRo07} or a wavelet basis~\cite{KrWa14}. The resulting methods cannot be considered to be fast transforms per se, as the computation time required for most reconstruction procedures scales at least linearly in the number of candidates $|\Gamma|\gg s$.

For any sublinear time algorithm to recover Fourier sparse signals,
significant subsampling compared to $|\Gamma|$ must be an integral part.
For the case of $d=1$ spatial dimension, a randomized approach that achieves runtime scaling quadratically in terms of the sparsity $s$ (up to logarithmic factors in $|\Gamma|$) was presented in \cite{GiGuInMuSt02}.
A deterministic, combinatorial algorithm with runtime complexity scaling quadratically in the sparsity~$s$ was presented in~\cite{Iwe10, Iwe13}. Improved randomized algorithms with runtime scaling linearly in the sparsity were introduced in~\cite{GiMuSt05, Iwe10, AlInKaPr12, AlInKaPr12b, Iwe13} and an improved deterministic algorithm in~\cite{LaWaCh13}. The latter one was enhanced once more in \cite{ChLaWa16} to accomodate for samples perturbed by Gaussian noise.

Many of the aforementioned algorithms have been generalized to the multi-dimensional case $d\geq 2$, and new algorithms have been developed. In Table~\ref{tab:complexities_comparison}, we give an overview on existing  multi-dimensional sparse FFT approaches and the ones presented in this paper. We compare sample and computational complexities as well as the general settings and types of sampling strategies.

\begin{table}[tb]
  \begin{small}
	\begin{tabular}{p{2.5cm}p{2.8cm}p{3.4cm}p{5.1cm}}
		\toprule
		method & sample\newline complexity & computational\newline complexity & comments \\
		\midrule
		\cite{ChChWa19}
		& $d\,s$\newline (w.h.p.)
		& $d\,s\,\log s$\newline (w.h.p.)
		& random signal setting\newline deterministic sampling strategy\newline success w.h.p. 
		\vspace*{.5em}\\
			\cite{ChChWa19b}
		& $d\,s\, \log N_\Gamma$\newline (w.h.p.)
		& $d\,s\,\log s \log N_\Gamma$\newline (w.h.p.)
		& random signal setting\newline deterministic sampling strategy\newline success w.h.p.
		\vspace*{.5em}\\
		\cite{InKa14}
		&$d^{\OO{d}}\,s\,\log N_\Gamma$
		&$d^{\OO{d}}\,s\,\log N_\Gamma $\newline $+\,d^3\,N_\Gamma^d\, \log^3 N_\Gamma$
		&
		arbitrary signal\newline random sampling strategy \newline success w.h.p.
		\vspace*{.5em}\\
        \cite{Mo17}
        & $d^3\,s^2\,N_\Gamma \log N_\Gamma$
        & $d^3\,s^2\,N_\Gamma^2 \, \log N_\Gamma$
        & arbitrary signal\newline deterministic sampling strategy \newline guaranteed success
		\vspace*{.5em}\\
		\cite[Theorem~7]{Iwe13}
		& $d^4\,s^2\,\log^4 (N_\Gamma)$
		& $d^4\,s^2\,\log^4 (N_\Gamma)$
		& arbitrary signal\newline deterministic sampling strategy \newline guaranteed success
		\vspace*{.5em}\\
		\cite[Corollary~4]{Iwe13}
		& $d^4\,s\,\log^4 (N_\Gamma)$
		& $d^4\,s\,\log^4 (N_\Gamma)$
		& arbitrary signal\newline random subsampling of deterministic strategy \cite[Theorem~7]{Iwe13} \newline success w.h.p.
		\vspace*{.5em}\\
		\cite{KaPoVo17}
		& $d\,s^2\,N_\Gamma\log^2s$\newline (w.h.p.) \vspace*{.3em} \newline $d\,s^3\,N_\Gamma\log^3s$ \newline (w.c.)\vspace*{.5em}
		& $d^2\,s^2\,N_\Gamma\log^3s$\newline (w.h.p.) \vspace*{.3em} \newline $d^2\,s^3\,N_\Gamma\log^4s$ \newline (w.c.) \vspace*{.5em}
		& arbitrary signal\newline random sampling strategy \newline success w.h.p. \\
		\cite{ChIwKr18,ChIwVo19}
		& $d^5\,s^3\,\log^4 (d\,s\,N_\Gamma)$
		& $d^6\,s^5\,\log^4 (d\,s\,N_\Gamma)$
		& arbitrary signal\newline random sampling strategy \newline success w.h.p.
		\vspace*{.5em}\\
		Theorem~\ref{thm:intro}
		& $s\,\log|\Gamma|$\newline $\lesssim d\,s\,\log N_\Gamma$
		& $(s\, \log s+d\,|\Gamma|)\log|\Gamma|$\newline  $\lesssim d^2\,N_\Gamma^d\log N_\Gamma$ \newline for $s\,\log s\lesssim d\,N_\Gamma^d$%
		& arbitrary signal\newline random sampling strategy \newline success w.h.p.\\
		\vspace*{.5em}\\
		Section~\ref{sec:sparseFFT}
		& $d\,s^2\,\log^2 (d\,s\,N_\Gamma)$\newline$+\,d\,s\,N_\Gamma\log(d\,s)$
		& $d^2\,s^2\,N_\Gamma\,\log^3(d\,s\,N_\Gamma)$ \newline  (w.h.p.) \vspace*{.3em}
		\newline $d^2\,s^3\,N_\Gamma\,\log^3(d\,s\,N_\Gamma)$ \newline  (w.c.)
		& arbitrary signal\newline random sampling strategy \newline success w.h.p.\\
		\bottomrule
	\end{tabular}
  \end{small}
	\caption{Overview on multi-dimensional sparse FFT approaches with $N_\Gamma$ as defined in~\eqref{eq:def_NI}.}
	\label{tab:complexities_comparison}
\end{table}

The deterministic approach in~\cite{LaWaCh13} and its variant for noisy samples \cite{ChLaWa16} have been extended to high spatial dimensions~$d$ in~\cite{ChChWa19} and~\cite{ChChWa19b}, respectively, with runtime complexity $\OO{d\,s\log s}$ and sample complexity
$\OO{d\,s}$ on average (both with an additional factor of $\log N_\Gamma$ for the variant addressing noisy samples), cf.\ Table~\ref{tab:complexities_comparison}. The method is highly efficient since the average case sample complexity is best possible, as a frequency set $I\subset\Z^d$ of cardinality $|I|=s$ has $d\,s$ many entries, and correspondingly, the average case runtime complexity is best possible up to a logarithmic factor.
However, the approach uses a random signal setting, i.e., the expectations and success probabilities are computed with respect to active frequencies collected in $I$ that are assumed to be distributed uniformly 
at random in a full (finite) tensor product grid as candidate set~$\Gamma$. Especially, smaller structured subsets, such as (subsets of) hyperbolic crosses, do not fit into this model assumption. In contrast to this, the other methods mentioned in Table~\ref{tab:complexities_comparison} use an ``arbitrary signal'' model, i.e., the active frequencies can be an arbitrary subset of a suitable candidate set~$\Gamma$.

In~\cite{Mo17}, the author presents a deterministic and noise robust high-dimensional sparse FFT approach. The
corresponding complexities scale merely quadratic with respect to the sparsity~$s$. 
However, the computational complexity of the Fourier transform $\OO{d^3\,s^2\,N_\Gamma^2 \, \log N_\Gamma}$ and the computational complexity of the construction of the sampling set $\OO{d^6\,s^2\,N_\Gamma\log^3N_\Gamma}$ indicates that the approach has advantages only for moderate expansions~$N_\Gamma$ and moderate dimensions~$d$.

The randomized approach in~\cite{AlInKaPr12b} has also been generalized to Fourier transforms in constant dimension \cite{InKa14}, but due to dimensional scaling of $d^{\OO{d}}$ in the runtime as well as sample complexity, cf.\ Table~\ref{tab:complexities_comparison}, the approach will not be feasible in higher spatial dimensions~$d$.

In contrast, the deterministic approach in~\cite{Iwe13} has been shown to generalize to high spatial dimensions~$d$ with a quadratic scaling of the complexities in the sparsity~$s$ and a polynomial scaling $d^4$ in the spatial dimension, cf.\ Table~\ref{tab:complexities_comparison}, but may exhibit limitations with respect to its numerical stability on large candidate sets $\Gamma$, in particular in high spatial dimensions~$d$. The reason for this issue is a transformation of the $d$-dimensional frequency domain $\Gamma$, which is assumed to be a tensor product grid, to a one-dimensional frequency domain $[0,\tilde{N})$, where $\tilde{N}\gtrsim|\Gamma|$ is necessary in order to obtain a unique mapping between both frequency grids.
Subsequently, one applies a one-dimensional sparse FFT approach to the new one-dimensional problem for typically huge $\tilde{N}$, which suffers from numerical issues due to the fact that highly oscillating basis functions -- with then possibly neighboring one-dimensional frequencies -- are difficult to distinguish when using only a few samples.
Remarkably, the locations of the samples and the algorithm itself are fully deterministic, and the method is guaranteed to successfully detect all active frequencies (on a machine with sufficiently high numerical precision).
The related randomized version in~\cite{Iwe13} uses a random subsampling and achieves a linear scaling of the complexities in~$s$, but also suffers from potential numerical issues for large~$d$.

The aforementioned approaches share the main characteristic that they are specifically designed for full (finite) tensor product grids in frequency domain as candidate sets~$\Gamma$.
In many applications, however, the function space under consideration motivates a significantly reduced set~$\Gamma$, such as a hyperbolic cross, or possibly an unstructured set~$\Gamma$. So only those Fourier $s$-sparse functions need to be considered with active frequencies in the candidate set~$\Gamma$.

High-dimensional sparse fast Fourier transforms for such scenarios have been designed, for example, using rank-$1$ lattices \cite{PoVo14,KaPoVo17}.
These methods use a dimension-incremental pairing technique, which can also be found in \cite{Zi79,PoTa12}. We describe the general ideas behind this technique in more detail later in this section.
The required runtime complexity of the algorithm in~\cite{KaPoVo17} is $\OO{d^2\,s^2\,N_\Gamma\log^3s}$ with high probability, where $N_\Gamma$ is the expansion of $\Gamma$, cf.~\eqref{eq:def_NI}, and the sample complexity
is $\OO{d\,s^2\,N_\Gamma\log^2s}$, cf.\ Table~\ref{tab:complexities_comparison}.
Due to the mild dependence on the spatial dimension~$d$ and the sparsity~$s$, this approach is well suited for high-dimensional problems.

Furthermore, a much more general approach with runtime complexity
$\OO{d^6\,s^5\,\log^4 (d\,s\,N_\Gamma)}$
and sample complexity
$\OO{d^5\,s^3\,\log^4 (d\,s\,N_\Gamma)}$
was presented in \cite{ChIwKr18,ChIwVo19}, which has also been applied to more general tensor product bases. It is based on a similar dimension-incremental pairing technique as in \cite{PoVo14,KaPoVo17} and internally uses ideas from compressed sensing.

One main contribution of this paper is the design of a sample-efficient sparse fast Fourier transform for Fourier sparse functions in high spatial dimensions~$d$ with frequencies in a given arbitrary candidate set $\Gamma\subset\Z^d$, $|\Gamma|<\infty$, where in particular $\Gamma$ may be much smaller than a full tensor product grid. The runtime complexity exhibits an additive joint dependence on~$s$ and $|\Gamma|$ up to logarithmic factors and -- most importantly -- a sample complexity that is linear in $s$ times a logarithmic factor in~$|\Gamma|$.
We stress on the fact that the algorithm succeeds with high probability, where the failure probability does not depend on the structure of the candidate set $\Gamma$ or the active frequency set $I$ of sparsity $s$. In particular, the presented estimates still hold true even if the frequency set $I$ is arbitrarily clustered.
More precisely, our main result reads as follows.
 
\begin{Theorem}\label{thm:intro} (Main Theorem).
	We consider the frequency sets $I\subset\Z^d$ and $\Gamma\subset\Z^d$, $|\Gamma|\ge8$, fulfilling
	\begin{equation}
	I\subset\Gamma\subset
	\bigtimes_{j=1}^d[a_j,b_j], \qquad
	\max_{j=1,\ldots,d} b_j-a_j\le 10.33\,|I|,\label{eq:subset_I_Gamma}
	\end{equation}
    and the parameter $\delta\in(0,1)$.
	Then, there exists a sampling strategy based on random rank-1 lattices using less than
	$$ 37\,|I|\,(\ln |\Gamma|-\,\ln\delta)$$
	sampling values of a multivariate trigonometric polynomial $p\in\Pi_I$ in order to identify all frequencies $\boldk\in I$ belonging to the non-zero Fourier coefficients $\hat{p}_\boldk\neq 0$ with probability at least $1-\delta$.
	This identification and the computation of the Fourier coefficients $\hat{p}_\boldk$ can be realized by Algorithm~\ref{alg:compute_p_from_Gamma} in less than
	$$C\,(|I|\log |I|+d\,|\Gamma|)\,(\ln |\Gamma|-\,\ln\delta)$$
	arithmetic operations with an absolute constant $C>0$.
\end{Theorem}
\begin{Remark}
	In Theorem~\ref{thm:intro} we can relax the subset property~\eqref{eq:subset_I_Gamma} to \begin{equation*}
	\left(I\cup\Gamma\right)
	\subset
	\bigtimes_{j=1}^d[a_j,b_j], \qquad
	\max_{j=1,\ldots,d} b_j-a_j\le 10.33\,|I|,
	\end{equation*}
	with the consequence that only the frequencies $\boldk\in\Gamma\cap I$ can be identified.
	Moreover, the constant $10.33$ in~\eqref{eq:subset_I_Gamma} can be adapted for specific applications, see also Remark~\ref{rem:parameter_c}.
\end{Remark}

In addition to the significant improvement for arbitrary sampling sets, Theorem \ref{thm:intro} also yields significant acceleration for sparse Fourier transforms on a full tensor grid, the setup studied in many of the works discussed above. Namely, in high spatial dimensions, our result can be combined with a so-called dimension-incremental approach, which is also a  key ingredient of many of the currently best known algorithms for sparse fast transforms.
 The idea of such approaches is to identify the multi-indices corresponding to the frequencies of the non-zero Fourier coefficients component by component, where in each step one works with the candidate set $\Gamma$ consisting of all indices that agree with the previously identified partial component vectors, see Section~\ref{sec:sparseFFT} for more details. 
 
 To identify the components of the active multi-indices, we will use 
 Theorem \ref{thm:intro}, which is based on random rank-1 lattices. In that sense our work follows a similar strategy as in~\cite{KaPoVo17} (and its predecessor \cite{PoVo14}), which also uses rank-$1$ lattices in the identification step of a dimension-incremental approach. The sparse transforms designed in these works, however, exhibit a larger runtime complexity due to a suboptimal
 scaling of the sample complexity required in each dimension-incremental step. 
   More precisely, they are based on spatial discretizations of full frequency candidate sets without
  taking advantage of the knowledge that only a few frequencies within these candidate sets are active. 
With Algorithm~\ref{alg:detect_I_from_Gamma} below at hand, in contrast, one can explicitly take this additional information into account to obtain a constructive sampling strategy with an improved sample complexity in combination with a comparable runtime complexity as in~\cite{KaPoVo17}. This is summarized in the following Theorem.

\begin{Theorem}\label{thm:dimin} (Dimension-incremental strategy).
	We consider the frequency sets $\Gamma\subset\Z^d$, $|\Gamma|<\infty$, and $I\subset \Gamma$. %
    Moreover, we define $N_\Gamma$ as in~\eqref{eq:def_NI}, e.g., $\Gamma\subset \big[\ceil{-N_\Gamma/2},\floor{N_\Gamma/2}\big]^d\cap\Z^d$.
	Then, there exists a sampling strategy based on random rank-1 lattices using 
	\[
      \OO{d\,|I|\max(|I|,N_\Gamma)\log^2\frac{d\,|I|\,N_\Gamma}{\delta}}
    \]
	sampling values that allows with probability at least $1-\delta$ to identify all multivariate trigonometric polynomials $p$ supported on $I$, that is, those $p$ with all its non-zero Fourier coefficients corresponding to index vectors in $I$ (in other words, trigonometric polynomials with $\hat{p}_\boldk= 0$ for all $\boldk \in \Z^d \setminus I$).
	This identification and the computation of the active Fourier coefficients~$\hat{p}_\boldk$, $\boldk\in I$, can be realized by a combination of Algorithm~\ref{algo:sfft_general} and Algorithm~\ref{alg:compute_p_from_Gamma_s}. This method has a computational complexity of
    \[\OO{d^2\,|I|^2\,N_\Gamma\,\log^3\frac{d\,|I|\,N_\Gamma}{\delta}}\]
    with probability at least $1-\delta$ as well as
    $\OO{d^2\,|I|^3\,N_\Gamma\,\log^3\frac{d\,|I|\,N_\Gamma}{\delta}}$
    in the worst case.
\end{Theorem}

\begin{sloppypar}
To compare this runtime and sample complexity to those of alternative dimension-incremental strategies, we again refer to Table~\ref{tab:complexities_comparison}.
In addition to the approach~\cite{KaPoVo17} that is based on rank-1 lattices, this includes \cite{ChIwKr18,ChIwVo19} which use compressed sensing techniques for the identification step.
These compressed sensing based approaches use the special structure of the candidate set~$\Gamma$: When restricting to each entry of the index set, one just encounters a regular grid in a much lower dimensional space, so one can efficiently apply compressed sensing techniques. The entries are then combined using an appropriate pairing technique, also based on compressed sensing ideas. 
This approach is much more general in that it allows for much larger classes of bounded orthonormal systems than just Fourier, but for that reason, it also does not take optimal advantage of the structural properties of the Fourier transform and the possible computational speedup by FFT techniques that our approach is able to exploit. This explains the inferior dimensional scaling of the runtime complexity.
\end{sloppypar}

 Besides these works, dimension-incremental methods have also been studied in combination with Prony's method~\cite{PoTa12}.
Also, there are further Prony-like techniques available in higher dimensions, cf.~\cite{CuHoKnLe20}. 
There is no explicit analysis of runtime and sample complexity available, however, which is why we cannot include this approach in Table~\ref{tab:complexities_comparison}.
 In addition, the feasibility of Prony's method in this context is severely limited to low sparsities due to
 stability problems arising even for sparsity levels in the order of a few hundreds.

\section{Reconstructing multivariate trigonometric polynomials from samples along random rank-1 lattices}
\label{sec:main_result}

\subsection{Sparse FFT via rank-1 lattices -- background and previous works}
\label{sec:main_result:background}

The frequency identification method we present in this work is based on a class of well established cubature methods in higher dimensions, so-called rank-1 lattice rules, which are special quasi Monte Carlo methods.
Such methods consider rank-1 lattices as sampling sets, that is, sets of the form
\begin{equation} \label{eq:lattice}
\Lambda(\boldz,M):=\{j\boldz/M\bmod \boldone\colon j=0,\ldots,M-1 \} \subset \T^d,
\end{equation}
where
$\boldz\in\Z^d$ and $M\in\N$ are called generating vector and lattice size of $\Lambda(\boldz,M)$, respectively. A candidate for an approximate integral of a function over a high-dimensional cube is then computed as the average of its samples along this lattice.

Naturally, it depends on the function class under consideration whether and to which accuracy this candidate approximates the true integral. In this paper, we will apply this approach for the identification of a sparse trigonometric polynomial
\begin{equation}\label{eq:trig_pol}
p(\boldx):=\sum_{\boldk\in I}\hat{p}_\boldk \,\e^{2\pi\ii\boldk\cdot\boldx},
\end{equation}
where $\boldk\cdot\boldx:=\sum_{j=1}^dk_jx_j$ denotes the usual inner product in $\R^d$.
Recall that we are interested in the case that the multivariate trigonometric polynomial is sparse, i.e., the index set $I\subset\Z^d$ -- the {\em frequency set} --  not only has finite cardinality, $|I|<\infty$, but is also small.

Note that the {\em Fourier coefficient} $\hat{p}_\boldk$ with frequency $\boldk$ can be computed
by evaluating the integral
\begin{equation}
\hat{p}_\boldk:=\int_{\boldx\in\T^d}p(\boldx) \,\e^{-2\pi\ii\boldk\cdot\boldx}\mathrm{d}\boldx.\label{eq:pFC}
\end{equation}
As it turns out, a rank-1 lattice rule applied to the integrand in~\eqref{eq:pFC} often returns the correct Fourier coefficient.
Indeed, the  rank-1 lattice rule with rank-1 lattice $\Lambda(\boldz,M)$ as in \eqref{eq:lattice} yields 
\begin{align}
\hat{p}_\boldk^{\Lambda(\boldz,M)}&:=\frac{1}{M}\sum_{j=0}^{M-1}p\left(\frac{j}{M}\boldz\right)\e^{-2\pi\ii\frac{j}{M}\boldk\cdot\boldz}
=\sum_{\boldh\in I}\hat{p}_\boldh\frac{1}{M}\sum_{j=0}^{M-1}\e^{2\pi\ii\frac{j}{M}(\boldh-\boldk)\cdot\boldz}\nonumber\\
&=\sum_{\boldh\in I}\hat{p}_\boldh\,\delta_0((\boldh-\boldk)\cdot\boldz \bmod M)\nonumber\\
&=\sum_{\boldh\in(\Lambda^\perp(\boldz,M))\cap (I-\boldk)}\hat{p}_{\boldh+\boldk},\label{eq:r1l_aliasing_formula}
\end{align}
where 
$\delta_0(0):=1$, $\delta_0(l):=0$ for $l\neq 0$, and 
$\Lambda^\perp(\boldz,M):=\{\boldh\in\Z^d\colon\boldh\cdot\boldz\equiv 0\imod{M}\}$
is the integer dual lattice of the rank-1 lattice $\Lambda(\boldz,M)$. As we will see, for a random lattice $\Lambda$, $\Lambda^\perp(\boldz,M)$ will typically intersect $I-\boldk$ only in $\boldh=0$ for $\boldk\in I$ or not at all for  $\boldk\notin I$, so $\eqref{eq:r1l_aliasing_formula}$ typically consists just of the single summand $\hat{p}_\boldk$ or even none, as desired.

When the intersection consists of more elements than just $0$, we speak of {\em aliasing}. For a specific frequency, as observed recently by one of the authors in~\cite{Kae17},  this happens only with small probability for random lattices even with small lattice sizes $M \asymp|I|$. Consequently, one has unique reconstruction of at least some of the original Fourier coefficients $\hat{p}_\boldk$ with a certain probability. Nevertheless each realization of such a random lattice will yield aliasing effects for some fraction of its Fourier coefficient. This means that our goal of correctly determining for all candidate frequencies $\boldh\in \Gamma$ whether they are active or not is typically impossible with a single realization. Even to just identify a superset of the active frequencies (with some false-positives allowed) requires a lattice of larger cardinality.
More precisely, the samples along a (single) rank-1 lattice rule do not allow for the reconstruction of all original Fourier coefficients $\{\hat{p}_\boldk\}_{\boldk\in I}$ unless the rank-1 lattice size $M$ is on the order of
$M\gtrsim |I|^2$ in the worst case \cite{kaemmererdiss}. In addition to this disadvantage, a corresponding generating vector~$\boldz$ needs to be constructed, e.g.\ using a component--by--component approach, which can lead to large computational costs, cf.~\cite[Page 3]{KaPoVo17} for a detailed discussion on this topic.

That is why in this paper, we work with multiple rank-$1$ lattices.
Building upon the work of \cite{Kae17}, we develop a strategy to identify the frequencies belonging to the non-zero Fourier coefficients of the polynomial $p$ within a frequency candidate set $\Gamma\supset I$ in Section~\ref{sec:main_result:new}. To this end, we present some notation and technical basics from \cite{Kae17} in the following.

As we use a random approach to generate the multiple rank-1 lattices, our method will not be guaranteed to be successful in all cases but with a chosen high probability. Naturally, the success rate will depend on the candidate set~$\Gamma$ to some extent. In particular, a candidate set~$\Gamma$ that is compact and of small size yields superior performance as compared to one which is extremely large or very wide spread. To describe the nature of the frequency candidate set~$\Gamma$, we will work with its cardinality $|\Gamma|$ and its expansion 
\begin{equation}
N_\Gamma:=\max_{j=1,\ldots, d}\left\{\max_{\zb k\in \Gamma}k_j-\min_{\zb l\in \Gamma}l_j\right\}.\label{eq:def_NI}
\end{equation}
Here, $N_\Gamma$ is the smallest edge length of a cube containing the frequency candidate set $\Gamma$. That is, there exists some $\boldk\in\Z^d$ such that $\Gamma \subset \{\boldk+\boldh\in\Z^d\colon\boldh\in[0,N_\Gamma]^d\}$.

If  $\boldk, \boldk'\in\Gamma$ agree up to multiples of $M$, $\boldk \equiv \boldk' \operatorname{mod} M$, sampling on a rank-one lattice of the type $\Lambda(\boldz,M)$ cannot distinguish $\boldk$ and $  \boldk'$. In general, one can at best identify the equivalence class mod $M$. For that reason, it is often  beneficial to represent each of these equivalence classes by an element in $\{0,\dots, M-1\}$, i.e., to consider
$$
\Gamma_{\!\!\bmod M}:=\{\zb h_{\zb k}:=\zb k-M\lfloor \tfrac{\zb k}{M} \rfloor\colon \zb k\in \Gamma\},
$$
where we assume $M\in\N$, and to aim to identify $I_{\!\!\bmod M} \cap \Gamma_{\!\!\bmod M}$ using the techniques introduced in this paper. Under the assumption that $|I_{\!\!\bmod M} \cap \Gamma_{\!\!\bmod M}|=|I \cap \Gamma|$ (that will always hold for $M$ large enough), this also allows for the identification of $|I\cap \Gamma|$.
We refer the reader to \cite[Lem.~2.3]{Kae17} for further details.
Working with the modified candidate set is preferred as for $M\le N_\Gamma$, dealing with
$\Gamma_{\!\!\bmod M}$ instead of $\Gamma$ leads to a smaller expansion of the candidate set~$\Gamma$ under consideration and
allows -- in specific situations -- for a significantly reduced number of sampling values required for the suggested approach,  Algorithm~\ref{alg:detect_I_from_Gamma} below, to identify the frequency support $I$
of the multivariate trigonometric polynomial $p$.

Furthermore, for technical reasons, it will be beneficial to use rank-1 lattice sizes $M$ that are prime numbers in addition to ensuring that $|\Gamma_{\!\!\bmod M}|=|\Gamma|$, i.e., we consider
\[M\in P^\Gamma:=\{ M'\in\N\colon M' \text{ prime with }|\Gamma_{\!\!\bmod M'}|=|\Gamma|\}.\]
This set contains at least all primes greater than $N_\Gamma$.

At this point, we stress on the fact that we will exploit the advantageous aliasing formula~\eqref{eq:r1l_aliasing_formula} in order to construct our algorithm.
In particular, the strategy is to randomly choose multiple suitable rank-1 lattices $\Lambda(\boldz_1,M_1), \ldots, \Lambda(\boldz_L,M_L)$ and
to exploit the fact that for a fixed frequency problematic aliasing effects occur with a
fixed probability for each of the rank-1 lattices under consideration. In order to control the probability 
that these aliasing effects occur in a reasonably large proportion of the $L$ rank-1 lattices, we
need to choose $L$ sufficiently large.
We call such a collection
$$
\Lambda=\Lambda(\boldz_1,M_1,\ldots,\boldz_L,M_L):=\bigcup_{\ell=1}^L\Lambda(\boldz_{\ell},M_{\ell})
$$
of rank-1 lattices $\Lambda(\boldz_1,M_1), \ldots, \Lambda(\boldz_L,M_L)$ a multiple rank-1 lattice configuration.

The separate consideration of the aliasing effects of each single rank-1 lattice $\Lambda(\boldz_\ell,M_\ell)$, $\ell=1,\ldots,L$, allows for the separate computation of possibly active $\hat{p}_\boldk^{\Lambda(\boldz,M)}$, $\boldk\in\Gamma$.
The structure of a specific rank-1 lattice $\Lambda(\boldz,M)$ will then provide an efficient algorithm that computes all (aliased)
Fourier coefficients $\hat{p}_\boldk^{\Lambda(\boldz,M)}$, $\boldk\in \Gamma$, at once requiring a computational complexity of $\OO{M\log M+d|\Gamma|}$. This is highly efficient compared to applications of matrix vector products, cf.~\cite{kaemmererdiss}.

\subsection{New approach and reconstruction guarantee}
\label{sec:main_result:new}

Given the role of the candidate set~$\Gamma$, it should not come as a surprise that one can only exactly identify the frequency support~$I$ of a multivariate trigonometric polynomial~$p$ provided that $I\subset\Gamma$ .
In practical scenarios, however, errors cannot always be avoided. 
For example when applying the methods developed in this section as frequency identification steps in a dimension-incremental
approach, cf.~Section~\ref{sec:sparseFFT} below, the candidate sets arise from previous estimation steps. So it may happen that this assumption is violated at some point. For that reason we do not explicitly require $I\subset \Gamma$ and 
show that the method will identify the frequency support $I\cap\Gamma$ of $p$ within $\Gamma$
with high probability for a suitable choice of parameters -- even when other frequencies are present. Note, however, that these parameters depend on $I\cup\Gamma$, so they 
 are not straight forward to determine for $I\not\subset\Gamma$ and unknown $I$.

The main tools for our analysis are the following generalizations of~\cite[Lemma~3.1 and Theorem~3.2]{Kae17}. 
\begin{Lemma}\label{lem:probab_alias_single_rand_r1l}
	Let $I\subset\Z^d$ and $\Gamma\subset\Z^d$ be frequency sets of finite cardinalities, $|I|<\infty$ and $|\Gamma|<\infty$.
	We fix a frequency $\boldk\in I\cup\Gamma$ and choose a prime number $M$
	such that $|(I\cup\Gamma)_{\!\!\mod M}|=|(I\cup\Gamma)|$.
	In addition, we choose a generating vector $\boldz\in [0,M-1]^d\cap\Z^d$ uniformly at random. Then, with probability not greater than $\frac{|I|}{M}$, the frequency $\boldk$ aliases to at least one frequency from $I\setminus\{\boldk\}$, i.e.,
	$$
	\mathbb{P}\left(\hat{p}_\boldk^{\Lambda(\boldz,M)}\neq \hat{p}_\boldk\right)\le \frac{|I|}{M},
	$$
	where $\hat{p}_\boldk=0$ for all $\boldk\in\Gamma\setminus I$, cf.~\eqref{eq:pFC}.
\end{Lemma}
\begin{proof}
	We start with the case $\boldk\in\Gamma\setminus I$, and we build the frequency set $\tilde{I}=I\cup\{\boldk\}$.
	Since we have $|(I\cup\Gamma)_{\!\!\mod M}|=|I\cup\Gamma|$ and $\tilde{I}\subset (I\cup\Gamma)$, it follows that $|\tilde{I}_{\mod M}|=|\tilde{I}|$, and we apply \cite[Lemma~3.1]{Kae17} with the frequency set $\tilde{I}$, which yields that $\boldk$ aliases to at least one frequency from $\tilde{I}\setminus\{\boldk\}=I$ with probability at most $\frac{|\tilde{I}|-1}{M}=\frac{|I|}{M}$. \newline
	In the case $\boldk\in I$, we have $|I_{\bmod M}|=|I|$ since $|(I\cup\Gamma)_{\mod M}|=|I\cup\Gamma|$. Applying \cite[Lemma~3.1]{Kae17} yields that $\boldk$ aliases to at least one other frequency from $I\setminus\{\boldk\}$ with probability at most $\frac{|I|-1}{M} < \frac{|I|}{M}$.
\end{proof}

\begin{Theorem}
	\label{thm:gen_mr1l_ref}
	Consider a frequency set $I\subset\Z^d$ and a set $\Gamma\subset\Z^d$ of frequency candidates, $|I|<\infty$ and $|\Gamma|<\infty$. In addition, we fix $\delta\in (0,1)$ and $\nu\in(0,1/2]$.
	Moreover, we determine the numbers
	$$
	\lambda \geq c\,|I| \qquad \text{ with } c> \frac{1}{\nu},
	$$
	$$
	L \ge \left\lceil \frac{c(c-2)}{(c\nu-1)^2\ln(c-1)} \, (\ln |\Gamma|-\,\ln\delta) \right\rceil.
	$$
	We choose $L$ not necessarily distinct lattice sizes $M_\ell\in P^{I\cup\Gamma}$, $M_\ell>\lambda$. For each $M_\ell$, we choose a generating vector $\boldz_\ell\in [0,M_\ell-1]^d\cap\Z^d$ uniformly at random.
	Then, the probability that the frequency $\boldk\in I\cup\Gamma$ aliases to a frequency from~$I\setminus\{\boldk\}$ for at least $\ceil{\nu L}\in[1,L]$ many of the rank-1 lattices $\Lambda(\boldz_\ell,M_\ell)$ is less than $\delta/|\Gamma|$, i.e.,
	$$
		\mathbb{P}\left(\left\{\Big|\Big\{\ell\in\{1,\ldots,L\}\colon \hat{p}_\boldk^{\Lambda(\boldz_\ell,M_\ell)}\neq \hat{p}_\boldk\Big\}\Big|\ge\nu L\right\}\right)\le \frac{\delta}{|\Gamma|}.
	$$
\end{Theorem}
\begin{proof}
	For the fixed frequency $\boldk\in I\cup\Gamma$, we define the random variables
\begin{equation}
	Y_\ell^\boldk :=
	\begin{cases}
	0 &\colon \text{ $\boldk$ does not alias to a frequency from $I\setminus\{\boldk\}$ using $\Lambda(\boldz_\ell,M_\ell)$}, \\
	1 &\colon \text{ $\boldk$ aliases to at least one frequency from $I\setminus\{\boldk\}$ using $\Lambda(\boldz_\ell,M_\ell)$}. \\
	\end{cases}\label{eq:def_Ylk}
\end{equation}
	We distinguish two different cases. First, we consider $\boldk\in I$ for $|I|=1$. Then we achieve $Y_\ell^\boldk=0$ for each $\ell\in\{1,\ldots,L\}$ and $\mathbb{P}\left\{ \sum_{\ell=1}^{L} Y_\ell^\boldk \geq \nu L \right\}=0<\frac{\delta}{|\Gamma|}$.
	Second, we consider the cases $\boldk\in I$ in conjunction with $|I|>1$ or $\boldk\in\Gamma\setminus I$.
	The random variables $Y_\ell^\boldk$, $\ell=1,\ldots,L$, are independent with uniformly bounded mean $\mu_\ell$, $0<\mu_\ell\leq \frac{|I|}{M_\ell} < \frac{|I|}{\lambda} \leq \frac{1}{c}<\nu$, where the upper bounds hold due to Lemma~\ref{lem:probab_alias_single_rand_r1l}. The strict inequality $0<\mu_\ell$ holds since we obtain the equality $Y_\ell^\boldk=1$ for the admissible choice $\boldz_\ell=(0,\ldots,0)^\top$, which then implies $\mathbb{P}\left\{ Y_\ell^\boldk = 1 \right\}>0$.
	\newline
	Using Hoeffding's inequality, cf.\ \cite[Thm.~1]{Hoeff63}, and $\mu<\nu\le 1/2$, we obtain
	\begin{align*}
	\mathbb{P}\left\{ \sum_{\ell=1}^{L} Y_\ell^\boldk \geq \ceil{\nu L} \right\}
	&=\mathbb{P}\left\{ \sum_{\ell=1}^{L} Y_\ell^\boldk \geq \nu L \right\}
	=
	\mathbb{P}\left\{ L^{-1} \sum_{\ell=1}^{L} Y_\ell^\boldk - \mu \geq \nu - \mu \right\} \\
	&\leq
	\mathbb{P}\left\{ L^{-1} \sum_{\ell=1}^{L} Y_\ell^\boldk - \mu \geq \nu - 1/c \right\}
	\leq
	\e^{-L(\nu - 1/c)^2\frac{1}{1-2\mu}\,\ln\frac{1-\mu}{\mu}}. \\
	\intertext{Since $\frac{1}{1-2\mu}\,\ln\frac{1-\mu}{\mu}$ is strictly decreasing, we continue}
	& < 
	\e^{-L(c\nu-1)^2\frac{\ln{(c-1)}}{c(c-2)}}
	\leq
	\e^{\ln\delta - \ln |\Gamma|}
	= \frac{\delta}{|\Gamma|},
	\end{align*}
	where Hoeffding's inequality holds since $0<\nu-1/c<\nu-\mu\le1/2-\mu<1-\mu$.
\end{proof}

The remaining part of this section works out and analyzes an algorithm that classifies the frequencies of a
given set of frequency candidates $\Gamma$ in two disjoint sets, see Algorithm~\ref{alg:detect_I_from_Gamma}.
Its main idea is to apply $L$ one-dimensional FFTs of lengths $M_1,\ldots,M_L$ in lines \ref{line:alg:detect_I_from_Gamma:FFT:begin}--\ref{line:alg:detect_I_from_Gamma:FFT:end} 
using the sampling values of the unknown sparse signal $p(\boldx)=\sum_{\boldk\in I}\hat{p}_\boldk \,\e^{2\pi\ii\boldk\cdot\boldx}$.
This requires $\OO{L\,M\log M}$ arithmetic operations, where $M:=\max\{M_\ell\colon \ell=1,\ldots,L\}$.
We rearrange the results $\hat{g}_h^{(\ell)}$ from line~\ref{line:alg:detect_I_from_Gamma:FFT:computation} using $h=\boldk\cdot\boldz_\ell\bmod{M_\ell}$
in order to determine $L$ (potentially) aliased Fourier coefficients $\hat{p}_{\boldk}^{(\ell)}:=\hat{g}^{(\ell)}_{\boldk\cdot\boldz_\ell\bmod{M_\ell}}$, $\ell=1,\ldots,L$, for each $\boldk\in\Gamma$, which requires $\OO{d\,L\,|\Gamma|}$ arithmetic operations.
Subsequently in line~\ref{alg:detect_I_from_Gamma_line_delta0}, one counts for fixed frequency $\boldk\in\Gamma$ how many of the coefficients~$\hat{p}_{\boldk}^{(\ell)}$, $\ell=1,\ldots,L$, are non-zero. If the fraction of the $\hat{p}_{\boldk}^{(\ell)}$ that are non-zero
is above a certain threshold, the frequency~$\boldk$ is kept, as it is likely to be part of the frequency support $I$ of the unknown trigonometric polynomial $p$.
So one collects this frequency in a set of detected frequencies $\tilde{I}$. If the fraction of the coefficients $\hat{p}_{\boldk}^{(\ell)}$ that are non-zero, $\ell=1,\ldots,L$, is below the threshold, one discards the frequency $\boldk$.
This classification has a computational complexity of $\OO{L\,|\Gamma|}$.
Altogether, this yields a total computational complexity of $\OO{L(M\log M+d\,|\Gamma|)}$ for Algorithm~\ref{alg:detect_I_from_Gamma}.

\begin{algorithm}[tb]
	\caption{Detecting the frequency set of a trigonometric polynomial~$p$.}\label{alg:detect_I_from_Gamma}
	\begin{tabular}{p{2cm}p{5.0cm}p{7.25cm}}
		Input: 	& $\Gamma\subset\Z^d$ 									& set of frequency candidates\\
		& $\Lambda:=\bigcup_{\ell=1}^L\Lambda(\boldz_\ell,M_\ell)$ 	& sampling nodes\\
		& $\left\{p(\boldx)\colon\boldx\in\Lambda\right\}$		& sampling values of $p\in\Pi_I$ \\
		& $\nu\in(0,1/2]$											& decision parameter
	\end{tabular}
	\begin{algorithmic}[1]
		\For{$\ell=1 \textnormal{ to }L$}\label{line:alg:detect_I_from_Gamma:FFT:begin}
		\State Compute $\hat{g}_h^{(\ell)}:=\sum_{j=0}^{M_\ell-1}p\left(\frac{j}{M_\ell}\boldz\right)\e^{-2\pi\mathrm{i}\,jh/M_\ell}$, $h=0,\ldots,M_\ell-1$, using a 1d FFT. \label{line:alg:detect_I_from_Gamma:FFT:computation}
		\EndFor\label{line:alg:detect_I_from_Gamma:FFT:end}
		\State Set $\tilde{I}:=\emptyset$.
		\ForEach{$\boldk\in\Gamma$}
		\If{$\sum_{\ell=1}^L\Big(1-\delta_0\big(\hat{g}^{(\ell)}_{\boldk\cdot\boldz_\ell\bmod{M_\ell}}\big)\Big)\ge\nu L$}\label{alg:detect_I_from_Gamma_line_delta0}
		\State Update $\tilde{I}:=\tilde{I}\cup\{\boldk\}$.
		\EndIf
		\EndFor
	\end{algorithmic}
	\begin{tabular}{p{2cm}p{5.0cm}p{7.25cm}}
		Output: & $\tilde{I}$								& detected frequency set\\
		\cmidrule{1-3}
		\multicolumn{3}{p{15.1cm}}{Computational complexity: $\OO{L(M\log M+d\,|\Gamma|)}$, where $M:=\max\{M_\ell\colon \ell=1,\ldots,L\}$}
	\end{tabular}
\end{algorithm}

\begin{Remark}
	In Algorithm \ref{alg:detect_I_from_Gamma}, we use a comparison to zero on line \ref{alg:detect_I_from_Gamma_line_delta0} implemented as Kronecker delta function $\delta_0$. Clearly,
	numerical implementations should take into account numerical inaccuracies and utilize a suitable approximation of this function.
\end{Remark}

The next two lemmas estimate the failure probabilities of the proposed detection method
and, subsequently, Corollary~\ref{cor:prob_bound_reco_gamma_I} reveals the connection
to Theorem~\ref{thm:gen_mr1l_ref}. In the following, we widely use the random variables $Y_\ell^\boldk$, that have already been defined in \eqref{eq:def_Ylk}.

\begin{Lemma}\label{lem:false_pos} \textnormal{(Probability estimate of false positive)}\newline
	For $\boldk\in \Gamma\setminus I$ and the computed Fourier coeffcients $\hat{p}_\boldk^{(\ell)}:=\hat{g}^{(\ell)}_{\boldk\cdot\boldz_\ell\bmod{M_\ell}}$ in line~\ref{line:alg:detect_I_from_Gamma:FFT:computation} of Algorithm~\ref{alg:detect_I_from_Gamma}, we observe
	\begin{align}
	\mathbb{P}\left(\left\{|\{\ell\in\{1,\ldots,L\}\colon \hat{p}_\boldk^{(\ell)}\neq 0\}|\ge\nu L\right\}\right)\le \mathbb{P}\left(\left\{\sum_{\ell=1}^LY_\ell^\boldk\ge \nu L\right\}\right).
	\label{eq:prob_est_false_pos}
	\end{align}
\end{Lemma}
\begin{proof}
	Since $\boldk\in \Gamma\setminus I$, $\hat{p}_\boldk^{(\ell)}=0$ holds in each case where $\boldk$ does not alias to any frequency $\boldh\in I$. Therefore, we necessarily need aliasing in order to obtain
	$\hat{p}_\boldk^{(\ell)}\neq 0$, which yields
	$$
	\left\{\{\ell\in\{1,\ldots,L\}\colon \hat{p}_\boldk^{(\ell)}\neq 0\}\ge\nu L\right\}\subset
	\left\{\parbox[c]{8cm}{$\boldk$ aliased for at least $\nu L$ of the rank-1 lattices $\Lambda(\boldz_\ell,M_\ell)$, $\ell=1,\ldots,L$, to at least one $\boldh\in I$}\right\}.
	$$
\end{proof}

\begin{Lemma}\label{lem:false_neg} \textnormal{(Probability estimate of false negative)}\newline
	For $\boldk\in I\cap\Gamma$ and the computed Fourier coefficients $\hat{p}_\boldk^{(\ell)}:=\hat{g}^{(\ell)}_{\boldk\cdot\boldz_\ell\bmod{M_\ell}}$ in line~\ref{line:alg:detect_I_from_Gamma:FFT:computation} of Algorithm~\ref{alg:detect_I_from_Gamma}, we observe
	\begin{align}
	\mathbb{P}\left(\left\{|\{\ell\in\{1,\ldots,L\}\colon \hat{p}_\boldk^{(\ell)}\neq 0\}|<\nu L\right\}\right)\le \mathbb{P}\left(\left\{\sum_{\ell=1}^LY_\ell^\boldk\ge L-\nu L\right\}\right).
	\label{eq:prob_est_false_neg}
	\end{align}
\end{Lemma}
\begin{proof}
	We have
	\begin{align*}
	\mathbb{P}\left(\left\{|\{\ell\in\{1,\ldots,L\}\colon \hat{p}_\boldk^{(\ell)}\neq 0\}|<\nu L\right\}\right)&=
	\mathbb{P}\left(\left\{|\{\ell\in\{1,\ldots,L\}\colon \hat{p}_\boldk^{(\ell)}= 0\}|\ge L-\nu L\right\}\right)
	\intertext{and the inclusion}
	\left\{|\{\ell\in\{1,\ldots,L\}\colon \hat{p}_\boldk^{(\ell)}= 0\}|\ge L-\nu L\right\}&\subset\left\{\parbox[c]{6cm}{$\boldk$ aliased for at least $L-\nu L$ of the rank-1 lattices $\Lambda(\boldz_\ell,M_\ell)$, $\ell=1,\ldots,L$, to at least one $\boldh\in I\setminus\{\boldk\}$}\right\},
	\intertext{which holds since aliasing is necessary in order to even obtain $\hat{p}_\boldk^{(\ell)}\neq \hat{p}_\boldk\neq 0$. Consequently, we achieve}
	\mathbb{P}\left(\left\{|\{\ell\in\{1,\ldots,L\}\colon \hat{p}_\boldk^{(\ell)}\neq 0\}|<\nu L\right\}\right)&\le 
	\mathbb{P}\left(\left\{\parbox[c]{6cm}{$\boldk$ aliased for at least $L-\nu L$ of the rank-1 lattices $\Lambda(\boldz_\ell,M_\ell)$, $\ell=1,\ldots,L$, to at least one $\boldh\in I\setminus\{\boldk\}$}\right\}\right)\\
	&=\mathbb{P}\left(\left\{\sum_{\ell=1}^LY_\ell^\boldk\ge L-\nu L\right\}\right).
	\end{align*}
\end{proof}

As a consequence of the last two lemmas, we need upper bounds on
$$
\mathbb{P}\left(\left\{\sum_{\ell=1}^LY_\ell^\boldk\ge \nu L\right\}\right)
\qquad\text{and}\qquad
\mathbb{P}\left(\left\{\sum_{\ell=1}^LY_\ell^\boldk\ge (1-\nu) L\right\}\right),
$$
which are themselves upper bounds on the failure probabilities of the classification of $\boldk\in\Gamma\setminus I$ and $\boldk\in I$ in line~\ref{alg:detect_I_from_Gamma_line_delta0} of Algorithm~\ref{alg:detect_I_from_Gamma}, respectively. In order to apply Theorem~\ref{thm:gen_mr1l_ref}, we need to choose $\nu\in(0,1/2]$, observing that
\begin{align}
\mathbb{P}\left(\left\{\sum_{\ell=1}^LY_\ell^\boldk\ge (1-\nu) L\right\}\right)\le\mathbb{P}\left(\left\{\sum_{\ell=1}^LY_\ell^\boldk\ge \nu L\right\}\right).\label{eq:bound_nu_one_half}
\end{align}
This restriction on $\nu$ allows for the application of Theorem~\ref{thm:gen_mr1l_ref} which leads to the
following statement about Algorithm \ref{alg:detect_I_from_Gamma}.

\begin{Corollary}\label{cor:prob_bound_reco_gamma_I}
	Let $I\subset\Z^d$ and $\Gamma\subset\Z^d$ with $|I|<\infty$ and $ |\Gamma|<\infty$ be given. Moreover,
	we choose the parameters $\nu\in(0,1/2]$ and $c>1/\nu$. Furthermore, we fix $\delta\in(0,1)$,
	\begin{align*}
	L &:=\ceil{\frac{c(c-2)}{(c\nu-1)^2 \; \ln(c-1)} \, (\ln |\Gamma|-\,\ln\delta)}, \text{ and} \\
	M&:=\min\left\{p\in P^{I\cup\Gamma}\colon p> c|I|\right\}.
	\end{align*}
	Subsequently, we randomly choose $\boldz_1,\ldots,\boldz_\ell\in[0,M-1]^d\cap\Z^d$.
	In addition, we assume that the sampling points
	$\{(\boldx,p(\boldx))\in\T^d\times\C\colon \boldx\in\Lambda(\zb z_\ell,M), \ell=1,\ldots,L\}$ of the multivariate trigonometric polynomial
	$p(\boldx)=\sum_{\boldk\in I}\hat{p}_\boldk \,\e^{2\pi\ii\boldk\cdot\boldx}$, $\hat{p}_\boldk\neq 0$ for each $\boldk\in I$, are given.
	Then, the probability that the output $\tilde{I}$ of Algorithm \ref{alg:detect_I_from_Gamma} does not equal $\Gamma\cap I$ is less than $\delta$.
\end{Corollary}
\begin{proof}
	Applying Theorem \ref{thm:gen_mr1l_ref}, Lemma~\ref{lem:false_pos}, Lemma~\ref{lem:false_neg}, and the union bound yields
	\begin{eqnarray*}
		\mathbb{P}\left(\tilde{I}\neq \Gamma\cap I\right)\!\!\!\!&=&\!\!\!\!
		\mathbb{P}\left(\bigcup_{\boldk\in\Gamma\setminus I}\left\{|\{\ell\colon \hat{p}_\boldk^{(\ell)}\neq 0\}|\ge\nu L\right\}
		\cup\bigcup_{\boldk\in I\cap\Gamma}\left\{|\{\ell\colon \hat{p}_\boldk^{(\ell)}\neq 0\}|<\nu L\right\}
		\right)\\
		&\overset{\eqref{eq:prob_est_false_pos}\;\&\;\eqref{eq:prob_est_false_neg}}{\le}&\sum_{\boldk\in\Gamma\setminus I}\mathbb{P}\left(\left\{\sum_{\ell=1}^LY_\ell^\boldk\ge \nu L\right\}\right)+\sum_{\boldk\in I\cap\Gamma}\mathbb{P}\left(\left\{\sum_{\ell=1}^LY_\ell^\boldk\ge (1-\nu) L\right\}\right)\\
		&\overset{\nu\le 1/2 \;\&\; \eqref{eq:bound_nu_one_half}}{\le}&
		\sum_{\boldk\in\Gamma}\mathbb{P}\left(\left\{\sum_{\ell=1}^LY_\ell^\boldk\ge \nu L\right\}\right)\overset{\text{Thm.~\ref{thm:gen_mr1l_ref}}}{<}|\Gamma|\frac{\delta}{|\Gamma|}=\delta.
	\end{eqnarray*}
\end{proof}

\begin{algorithm}[tb]
	\caption{Detecting the frequency set and computing the Fourier coefficients of a trigonometric polynomial~$p$. (Algorithm~\ref{alg:detect_I_from_Gamma} with $\nu:=1/2$ and computation of Fourier coefficients).}\label{alg:compute_p_from_Gamma}
	\begin{tabular}{p{2cm}p{5.0cm}p{7.25cm}}
		Input: 	& $\Gamma\subset\Z^d$ 									& set of frequency candidates\\
		& $\Lambda:=\bigcup_{\ell=1}^L\Lambda(\boldz_\ell,M_\ell)$ 	& sampling nodes, $L$ odd\\
		& $\left\{p(\boldx)\colon\boldx\in\Lambda\right\}$		& sampling values of $p\in\Pi_I$ \\
	\end{tabular}
	\begin{algorithmic}[1]
		\For{$\ell=1 \textnormal{ to }L$}
		\State Compute $\hat{g}_h^{(\ell)}:=\frac{1}{M_\ell}\sum_{j=0}^{M_\ell-1}p\left(\frac{j}{M_\ell}\boldz\right)\e^{-2\pi\mathrm{i}\,jh/M_\ell}$, $h=0,\ldots,M_\ell-1$, using a 1d FFT.\label{alg:2_line_hatg}
		\EndFor
		\State Set $\tilde{I}:=\emptyset$.
		\ForEach{$\boldk\in\Gamma$}
		\If{$\sum_{\ell=1}^L\Big(1-\delta_0\big(\hat{g}^{(\ell)}_{\boldk\cdot\boldz_\ell\bmod{M_\ell}}\big)\Big)\ge L/2$}
		\State Update $\tilde{I}:=\tilde{I}\cup\{\boldk\}$.
        \State \parbox[t]{\linewidth-\algorithmicindent-\algorithmicindent}{With $\hat{p}_\boldk^{(\ell)}:=\hat{g}^{(\ell)}_{\boldk\cdot\boldz_\ell\bmod{M_\ell}}$, compute  \newline   $\check{p}_\boldk:=\operatorname{median}\left\{\operatorname{Re}\big(\hat{p}_\boldk^{(\ell)}\big)\colon\ell=1,\ldots,L\right\} + \mathrm{i} \cdot \operatorname{median}\left\{\operatorname{Im}\big(\hat{p}_\boldk^{(\ell)}\big)\colon\ell=1,\ldots,L\right\}$. \label{line:alg:compute_p_from_Gamma:median}}
		\EndIf
		\EndFor
	\end{algorithmic}
	\begin{tabular}{p{2cm}p{5.0cm}p{7.25cm}}
		Output: & $\tilde{I}$								& detected frequency set\\
		& $(\check{p}_\boldk)_{\boldk\in\tilde{I}}$	& Fourier coefficients belonging to $\tilde{I}$ %
		\\    \cmidrule{1-3}
		\multicolumn{3}{p{15.1cm}}{Computational complexity: $\OO{L(M\log M+d\,|\Gamma|)}$, where $M=\max\{M_\ell\colon \ell=1,\ldots,L\}$}
	\end{tabular}
\end{algorithm}

Up to now, we considered only the classification of the non-zero frequencies.
Specific parameter choices allow for the additional computation of the unknown Fourier coefficients.
To this end, we compute
\begin{align}
\check{p}_\boldk:=\operatorname{median}\left\{\operatorname{Re}\big(\hat{p}_\boldk^{(\ell)}\big)\colon\ell=1,\ldots,L\right\} + \mathrm{i} \cdot \operatorname{median}\left\{\operatorname{Im}\big(\hat{p}_\boldk^{(\ell)}\big)\colon\ell=1,\ldots,L\right\},
\label{eq:def_median_check_p}
\end{align}
where $\hat{p}_\boldk^{(\ell)}:=\hat{g}^{(\ell)}_{\boldk\cdot\boldz_\ell\bmod{M_\ell}}$.
Choosing $\nu=1/2$ and $L$ odd, we obtain
$$\check{p}_\boldk=
\begin{cases}
0 &\colon {\boldk\in\Gamma\setminus I},\\
\hat{p}_\boldk &\colon {\boldk\in I\cap\Gamma},
\end{cases}
$$
with probability at least $\eta:=1-\mathbb{P}\left(\left\{\sum_{\ell=1}^LY_\ell^\boldk\ge L/2\right\}\right)$ since there
is no aliasing for at least $L/2$, i.e., $\frac{L+1}{2}$, rank-1 lattices on the frequency $\boldk$ with at least this probability.
Accordingly, with probability $\eta$, we determine for at least $\frac{L+1}{2}$ different $\ell$ the correct values $\hat{p}_\boldk^{(\ell)}=\hat{p}_\boldk$, which implies the median is exactly this value.
Algorithm~\ref{alg:compute_p_from_Gamma} presents the resulting strategy for detecting the active frequencies $\boldk\in\Gamma$ as well as
computing all medians of the sets $\{\hat{p}_\boldk^{(\ell)}\colon \ell=1,\ldots, L\}$, $\boldk\in\tilde{I}$,
as corresponding Fourier coefficients.

The differences to Algorithm~\ref{alg:detect_I_from_Gamma} are the fixed parameter choice $\nu:=1/2$ and the additional line~\ref{line:alg:compute_p_from_Gamma:median} in Algorithm~\ref{alg:compute_p_from_Gamma}, which computes the Fourier coefficients $\check{p}_\boldk$, $\boldk\in\tilde{I}$, as medians using the method from~\cite{BlFlPrRiTa73}.
Since $\tilde{I}\subset\Gamma$, Algorithm~\ref{alg:compute_p_from_Gamma} requires at most $\OO{L\,|\Gamma|}$ additional arithmetic operations. Thus, we obtain the same computational complexities for Algorithms~\ref{alg:detect_I_from_Gamma} and~\ref{alg:compute_p_from_Gamma}.

\begin{Corollary}\label{cor:prob_bound_reco_trig_poly}
	Let $I\subset\Z^d$ and $\Gamma\subset\Z^d$ with $|I|<\infty$ and $ |\Gamma|<\infty$ be given. Moreover,
	we choose $\nu:=1/2$ and $c>2$. In addition, we fix $\delta\in(0,1)$,
	\begin{align*}
	L &:=\min\left\{n\in 2\N+1\colon n\ge \frac{4c}{(c-2)\,\ln(c-1)} \, (\ln |\Gamma|-\,\ln\delta)\right\}, \text{ and} \\
	M&:=
	\min\left\{p\in P^{I\cup\Gamma}\colon p> c|I|\right\}.
	\end{align*}
	Subsequently, we randomly choose $\boldz_1,\ldots,\boldz_\ell\in[0,M-1]^d\cap\Z^d$.
	In addition, we assume that the sampling points
	$\{(\boldx,p(\boldx))\in\T^d\times\C\colon \boldx\in\Lambda(\zb z_\ell,M%
	), \ell=1,\ldots,L\}$ of the multivariate trigonometric polynomial
	$p(\boldx)=\sum_{\boldk\in I}\hat{p}_\boldk\e^{2\pi\ii\boldk\cdot\boldx}$, $\hat{p}_\boldk\neq 0$ for each $\boldk\in I$, are given.
	Then, the probability that the output $(\check{p}_\boldk)_{\boldk\in\tilde{I}}$ of Algorithm \ref{alg:compute_p_from_Gamma} does not match the correct Fourier coefficients, i.e., $\tilde{I}\neq I\cap\Gamma$ or $\check{p}_\boldk\neq\hat{p}_\boldk$ for at least one $\boldk\in \tilde{I}$,
	is bounded from above by $\delta$.
\end{Corollary}
\begin{proof}
	We follow the proof of Corollary \ref{cor:prob_bound_reco_gamma_I} for estimating the probability, and we take into account the argumentations on the medians $\check{p}_\boldk$ immediately above this Corollary.
\end{proof}

\begin{Remark}\label{rem:parameter_c}
	The parameters in Corollary~\ref{cor:prob_bound_reco_trig_poly} need to be chosen carefully depending on the specific application. In particular, it might be useful to increase the parameter~$c$, which leads to possibly larger lattice sizes $M$ but smaller numbers $L$ of used rank-1 lattices.
\end{Remark}

Theorem~\ref{thm:intro} states the result of Corollary~\ref{cor:prob_bound_reco_trig_poly} for the specific parameter choice $c:=10.33$ and the restriction on the frequency sets given in~\eqref{eq:subset_I_Gamma}. We now give its proof.

\begin{proof}[Proof of Theorem.~\ref{thm:intro}]
	We apply Corollary~\ref{cor:prob_bound_reco_trig_poly} setting the parameter $c:=10.33$.  
	Since the subset property \eqref{eq:subset_I_Gamma} of the frequency sets is assumed,
	$P^{I\cup\Gamma}$ contains each prime number $p$ larger than $c|I|$.\newline
	The ratio $\frac{M}{c|I|}$ of the smallest prime number $M$ larger than $c|I|$ is bounded from above
	by $\frac{M}{c|I|}\le\frac{127}{11\cdot10.33}$, which can be proven using \cite[Prop.~5.4]{Du18}
	and the	calculation of finitely many instances of this ratio.
	Accordingly, the lattice sizes are bounded from above by $M\le \frac{127}{11}|I|$.
	In addition, for $|\Gamma|\geq 8>e^2$, the number of lattices $L$ is bounded from above by
	\begin{align*}
	L&\le \frac{4c}{(c-2)\ln\,(c-1)} \, (\ln |\Gamma|-\,\ln\delta)+2\\
	&\le \frac{4c}{(c-2)\,\ln(c-1)}  \left(\ln |\Gamma|-\,\ln\delta\right)+ \frac{2}{\ln{8}}(\ln|\Gamma|-\ln\delta)\\
	&\le 3.183\left(\ln |\Gamma|-\,\ln\delta\right)
	\end{align*}
	where in the last inequality we used that $\delta<1$. This yields a total number of samples $M\,L<37|I|\left(\ln |\Gamma|-\,\ln\delta\right)$.
	According to Corollary~\ref{cor:prob_bound_reco_trig_poly}, we apply Algorithm~\ref{alg:compute_p_from_Gamma} and obtain the non-zero ones of the Fourier coefficients $\hat{p}_\boldk$, $\boldk\in I=I\cap\Gamma$, with probability at least $1-\delta$.
\end{proof}

\begin{Remark}\label{rem:use_r1l_info_or_ls}

Clearly, the reconstruction guarantees in Theorem~\ref{thm:intro} as well as Corollary~\ref{cor:prob_bound_reco_trig_poly}
hold with probability of at least $1-\delta$. However, from a practical point of view, it is reasonable to expect that the output $\tilde{I}$ of Algorithms~\ref{alg:detect_I_from_Gamma} and~\ref{alg:compute_p_from_Gamma}
contains $I$, i.e., $I\subset\tilde{I}$, with significantly higher probability. 
On the one hand, the estimate \eqref{eq:prob_est_false_neg} in Lemma~\ref{lem:false_neg} is very rough since we
estimate the probability of cancellations of aliasing Fourier coefficients just by the probability of aliasing frequencies, which are widely differing events in fact.
On the other hand, the probability that $I\not\subset\tilde{I}$ holds, i.e., that we observe false negatives, is bounded from above by $\frac{|I|}{|\Gamma|}\delta$, cf.\ the proof of Corollary~\ref{cor:prob_bound_reco_gamma_I}.
Moreover, the probability that we will observe false positives seems to be significantly higher, cf.\ Section~\ref{ssec:num1}, which also fits well to the theoretical estimate on the failure probability in Corollary~\ref{cor:prob_bound_reco_gamma_I} since the upper bound on the probability of observing false positives
contributes much more to the upper bound $\delta$ of the total failure probability when one assumes that $|I|\ll|\Gamma|$.

From that point of view, the probability of $I\subset\tilde{I}$ in conjunction with, for instance, $|\tilde{I}|\le 2\,|I|$ may be considerably larger than the probability we estimated in Theorem~\ref{thm:intro}.
Even in these cases, one has a reasonable chance to entirely identify the trigonometric polynomial $p$
from the already known sampling values.

An improvement strategy is to postprocess the frequency set $\tilde{I}$ and Fourier coefficients~$\check{p}_\boldk$, $\boldk\in\tilde{I}$, obtained as output of Algorithm~\ref{alg:compute_p_from_Gamma}. One considers $\tilde{I}$ as the frequency
set of the polynomial $p$ and one computes the corresponding Fourier coefficients~$\hat{p}_{\boldk}$, $\boldk\in\tilde{I}$, from the already available sampling values
along the fixed rank-1 lattices $\Lambda(\boldz_1,M_1),\ldots,\Lambda(\boldz_L,M_L)$ using \cite[Algorithm~2]{KaVo19} or, even more general, one applies a least squares method. If this set of rank-1 lattices provides a spatial discretization of the space $\Pi_{\tilde{I}}$ of trigonometric polynomials, cf.~\cite{Kae17} for details, and $\tilde{I}\supset I$ holds, we will entirely identify the trigonometric polynomial~$p$. Otherwise, we can use the following hybrid approach. For each rank-1 lattice $\Lambda(\boldz_\ell,M_\ell)$, we determine the frequency set 
$
\tilde{I}^{(\ell)}:=\{\boldk\in\tilde{I}\colon \boldk\cdot\boldz_\ell \not\equiv \boldh\cdot\boldz_\ell \imod{M_\ell} \forall\boldh\in\tilde{I}\setminus\{\boldk\}\}
$
belonging to reconstructable Fourier coefficients $\hat{p}_\boldk$ for any $p\in\Pi_{\tilde{I}}$ and set $\check{p}_\boldk:=\hat{p}_\boldk^{\Lambda(\boldz_\ell,M_\ell)}$ for $\boldk\in\tilde{I}^{(\ell)}$. Whenever a frequency $\boldk\in\tilde{I}$ is contained in more than one $\tilde{I}^{(\ell)}$, we set the Fourier coefficient~$\check{p}_\boldk$ to the average of the corresponding coefficients~$\hat{p}_\boldk^{\Lambda(\boldz_\ell,M_\ell)}$. If there should be frequencies $\boldk\in\tilde{I}\setminus\big(\cup_{\ell=1}^L\tilde{I}^{(\ell)}\big)$, then we just use the Fourier coefficients~$\check{p}_\boldk$ computed in line~\ref{line:alg:compute_p_from_Gamma:median} of Algorithm~\ref{alg:compute_p_from_Gamma} as a fall-back.

Numerical tests, cf.\ Section~\ref{ssec:num1}, confirm that this postprocessing strategy works very well. However,
we cannot directly apply the theoretical results from~\cite{Kae17}, since the frequency set~$\tilde{I}$ already depends on the aliasing effects of the rank-1 lattices $\Lambda(\boldz_1,M_1),\ldots,\Lambda(\boldz_L,M_L)$.
\end{Remark}

\section{A sparse FFT for high-dimensional data}\label{sec:sparseFFT}

\subsection{Dimension-incremental sparse FFT -- background and previous works}

As already mentioned in the introduction, one of the main motivations for our algorithm is to apply it as a key ingredient  of a dimension-incremental algorithm for the high-dimensional fast Fourier transform.  
Such algorithms have received considerable attention in recent years 
due to their excellent and reliable applicability in high-dimensional settings, cf.,\ e.g.,\ \cite{PoVo14, KaPoVo17, ChIwKr18, ChIwVo19}.

To precisely formulate the algorithmic framework of these approaches that will also form the basis of our sparse FFT, 
we recall some notation:
We denote the  multivariate trigonometric polynomial that we aim to recover by
\begin{equation}
p\colon \T^d\rightarrow\C, \qquad p(\boldx) := \sum_{\boldk\in I} \hat{p}_\boldk \, \e^{2\pi\ii\boldk\cdot\boldx}, \quad \hat{p}_\boldk\in\C\setminus\{z\in\C\colon |z|<3 \theta\},
\label{eq:trig_poly_p}
\end{equation}
for some frequency set $I=\operatorname{supp}\,\hat{p}\subset\Z^d$, $s:=|I|<\infty$.
Furthermore, we assume that we know a (possibly very large) candidate set~$\Gamma\subset\Z^d$ of finite cardinality for~$I$, i.e., $I\subset\Gamma$ and $|\Gamma|<\infty$.

In this notation our key dimensional-incremental strategy, in analogy to a number of recent approaches in the literature, such as \cite{ChIwKr18,ChIwVo19}, amounts to first identifying $I^{(t)}$ as a superset of the projection
$$\mathcal{P}_t (I):=\{k_t\in\Z\colon \boldk=\left(k_1,\ldots,k_d\right)^\top\in I\}$$
to the $t$th component of the active frequencies in $I$
for each coordinate $t\in\{1, \dots, d\}$, and then incrementally combining ("pairing") the elements of the different $I^{(t)}$ to iteratively obtain a superset $I^{(1,\ldots,t)}$ of $$\mathcal{P}_{1,\ldots,t} (I):=\{(k_1,\ldots,k_t)\in\Z^t\colon\boldk=(k_1,\ldots,k_d)^\top\in I\},$$ which will eventually yield a set of multi-indices $I^{(1,\ldots,d)} \supset \mathcal{P}_{1,\ldots,d} (I)=I$ that contains the active ones.

In each of these pairing steps, one
determines the
candidates for the indices in $I^{(1,\ldots,t)}$ by appropriately combining the (already identified) frequency set $I^{(1,\ldots,t-1)}$ for $\mathcal{P}_{1,\ldots,t-1} (I)$ with the (already identified) frequency set $I^{(t)}$ for $\mathcal{P}_t (I)$. That is, one aims to find  $|I|$ active elements from the larger candidate set $(I^{(1,\ldots,t-1)}\times I^{(t)})\cap\mathcal{P}_{1,\ldots,t}(\Gamma) =:J_t$. This problem is exactly of the form studied in the previous section, and basically this step is where the various approaches that have been proposed in the literature differ. Treating this step as a black box, say Algorithm $\operatorname{A}$, we can formulate the dimension-incremental approach to high-dimensional fast transforms as a meta-algorithm, which is summarized in Algorithm~\ref{algo:sfft_general}.

To apply our algorithm for a wider class of scenarios, we introduced the two additional parameters $\theta$ and $s_\mathrm{local}$ for enhanced stability and robustness. 
The threshold parameter $\theta$ is meant to account for scenarios where due to noise or model mismatch, some linear combination of non-significant coefficients is not evaluated to exactly zero, but just approximately zero.
More precisely, one is not aiming to identify all non-zero coefficients, but rather only those with absolute values of at least~$\theta$. Naturally the parameter $\theta$ should be chosen smaller than the values to be expected, see also \eqref{eq:trig_poly_p}.

 The parameter $s_\mathrm{local}$ is related to the observation that allowing for a moderate number of false positives can significantly improve the success rate for recovery in practice. If Step~2c of Algorithm~\ref{algo:sfft_general} is designed to identify only exactly $s$ coefficients, then a false positive 
can lead to the situation where an active frequency no longer corresponds to any of the $s$ largest coefficients, and hence is not identified.

For that reason it is useful to collect $s_\mathrm{local}>s$ frequencies in the intermediate steps of Algorithm~\ref{algo:sfft_general} in practice, %
and we incorporate this feature in Algorithm~\ref{algo:sfft_general} by using  $s_\mathrm{local}:=2s$ as the default choice. 

\begin{algorithm}[p]
	\caption{Dimension-incremental reconstruction of a multivariate trigonometric polynomial $p$ from sampling values.}\label{algo:sfft_general}
  \begin{small}
	\begin{tabular}{p{1.2cm}p{2.1cm}p{11cm}}
		Input:  %
		& $\Gamma\subset\Z^d$ \hfill & search space in frequency domain, candidate set for $I=\mathrm{supp}\,\hat{p}$ \\
		& $p(\circ)$ & trigonometric polynomial $p$ as black box (function handle) \\
		& $s,s_\mathrm{local}\in\N$ & sparsity parameter, $s\leq s_\mathrm{local}$  ($s_\mathrm{local}:=2s$ by default) \\  
		& Algorithm $\operatorname{A}$ & {\bfseries efficient algorithm} $\operatorname{A}$ that guarantees the identification of the frequency support of each $s_\mathrm{local}$-sparse trigonometric polynomial in $\Pi_{J_t}$
        w.h.p., cf.\ Section~\ref{sec:complexity}, and computes the Fourier coeffcients \\
		& $\theta\in\R^+$ & absolute threshold (e.g. close to machine precision) \\
		& $r\in\N$ & number of detection iterations%
	\end{tabular}
	\vspace{-0.6em}
	\begin{algorithmic}
		\item[(Step 1)] [Single frequency component identification]\label{algo:sfft_general:step1}
		\For{$t:=1,\ldots,d$}
		\State Set $K_t:=\max(\mathcal{P}_t(\Gamma))-\min(\mathcal{P}_t(\Gamma))+1$, $I^{(t)}:=\emptyset$.
		\For{$i:=1,\ldots,r$}
		\State Choose $x_j'\in\T$, $j\in\{1,\ldots,d\}\setminus\{t\}$ uniformly at random.
		\State Set $\boldx^{(\ell)}:=\left(x_1^{(\ell)},\ldots, x_d^{(\ell)}\right)^\top$, $x_j^{(\ell)}:=
		\begin{cases}
		\ell/K_t,& j=t,\\
		x_j',& j\neq t,
		\end{cases}$ \quad
		for all $\ell=0,\ldots,K_t-1$.
		\State \vspace{0.2em}
		\State Compute $\tilde{\hat{p}}_{t,k_t}:=\frac{1}{K_t} \sum_{\ell=0}^{K_t-1} p\left(\boldx^{(\ell)}\right) \,\mathrm{e}^{-2\pi\mathrm{i} \ell k_t / K_t}$, $k_t\in\mathcal{P}_t(\Gamma)$, via FFT.
		\State Set $I^{(t)}:=I^{(t)}\cup\{k_t\in\mathcal{P}_t(\Gamma)\colon \tilde{\hat{p}}_{t,k_t} \text{ is among the largest } s_\mathrm{local} \text{ (in absolute value)}$
		\State \qquad\qquad\qquad\qquad $\text{ elements of } \{\tilde{\hat{p}}_{t,j}\}_{j\in\mathcal{P}_t(\Gamma)} \text{ and } \vert\tilde{\hat{p}}_{t,k_t}\vert \geq \theta
		\}$.   \EndFor{\, $i$}
		\EndFor{\, $t$}
		\item[(Step 2)] [Coupling frequency components identification]
		\For{$t:=2,\ldots,d$}
		\State If $t<d$, set $\tilde r:=r$ and $\tilde s:=s_\mathrm{local}$, otherwise $\tilde r:=1$ and $\tilde s:=s$.
		Set $I^{(1,\ldots,t)}:=\emptyset$.
		\For{$i:=1,\ldots,\tilde r$}
		\State Choose components $x_{t+1}',\ldots,x_d'\in\T$ of sampling nodes uniformly at random.
		\item[(Step 2a)]
		\State \parbox[t]{\linewidth-\algorithmicindent-\algorithmicindent}{Generate a sampling set $\mathcal{X}\subset\T^t$ for 
			$ J_t:=(I^{(1,\ldots,t-1)} \times I^{(t)})\cap\mathcal{P}_{(1,\ldots,t)}(\Gamma)$ that allows for the application of  Algorithm $\operatorname{A}$.
			Set $\mathcal{X}_{t,i}:=\{\boldx:=(\boldsymbol{\tilde{x}},x_{t+1}',\ldots,x_d')\colon \boldsymbol{\tilde{x}}\in\mathcal{X} \}\subset\T^d$.}
		\item[(Step 2b)]
		\State \parbox[t]{\linewidth-\algorithmicindent-\algorithmicindent}{Sample $p$ along the nodes of the sampling set $\mathcal{X}_{t,i}$.}
		\item[(Step 2c)]
		\State \parbox[t]{\linewidth-\algorithmicindent-\algorithmicindent}{Apply Algorithm~$\operatorname{A}$ to obtain the support $\tilde{J}_{t,i}\subset J_t$, $|\tilde{J}_{t,i}|\le \tilde{s}$, of frequencies belonging to the at most $\tilde{s}$ largest Fourier coefficients, each larger than $\theta$ in absolute value, using the sampling values~$p(\boldx_j)$, $\boldx_j\in\mathcal{X}_{t,i}$.}
		\item[(Step 2d)]
		\State Set $I^{(1,\ldots,t)}:=I^{(1,\ldots,t)}\cup \tilde{J}_{t,i}$.
		\EndFor{\, $i$}
		\EndFor{\, $t$}
		\item[(Step 3)] [Computation of Fourier coefficients]
		\State Generate a sampling set $\mathcal{X}\subset\T^d$ for $I^{(1,\ldots,d)}$ such that the corresponding Fourier matrix $\boldsymbol{A}(\mathcal{X},I^{(1,\ldots,d)})$ is of full column rank and its pseudoinverse can be applied {\bfseries efficiently}.
		Compute the corresponding Fourier coefficients $\left(\tilde{\hat{p}}_{(1,\ldots,d),\boldk}\right)_{\boldk\in I^{(1,\ldots,d)}}$.
		\State Set $\tilde{I}:=\Big\{\boldk\in I^{(1,\ldots,d)}\colon \tilde{\hat{p}}_{(1,\ldots,d),\boldk} \text{ is among the largest } s \text{ (in absolute value) elements of }$
		\State \qquad\qquad\qquad\qquad\qquad $\{\tilde{\hat{p}}_{(1,\ldots,d),\boldj}\}_{\boldj\in I^{(1,\ldots,d)}} \text{ and } |\tilde{\hat{p}}_{(1,\ldots,d),\boldk}| \geq \theta \Big\}$, $\boldsymbol{\tilde{\hat{p}}}:=\left(\tilde{\hat{p}}_{(1,\ldots,d),\boldk}\right)_{\boldk\in \tilde{I}}$.
	\end{algorithmic}

	\begin{tabular}{p{1.2cm}p{2.0cm}p{11.1cm}}
		Output: & $\tilde{I}\subset\Gamma\subset\Z^d$ & index set of detected frequencies, $|\tilde{I}|\leq \min \{s,|\Gamma|\}$ \\
		& $\boldsymbol{\tilde{\hat{p}}}\in\C^{\vert \tilde{I}\vert}$ & corresponding Fourier coefficients, $|\tilde{\hat{p}}_{(1,\ldots,d),\boldk}| \geq \theta$ \\
	\end{tabular}
  \end{small}
\end{algorithm}

In general, a major computational bottleneck of such dimension-incremental approaches is the iterated signal reconstruction from sampling values, (i.e., Step 2 in the formalization of Algorithm~\ref{algo:sfft_general}). In the simplest case, each such instance can be thought of as computing a pseudoinverse of a matrix. Since the matrices to be inverted change in every step of the dimension-incremental algorithm, a procedure for highly efficient computation of such a pseudoinverse is the key to the computational efficiency of the entire method. If one has FFT-like algorithms at one's disposal, the computational complexity will be tremendously reduced as compared to algorithms entirely based on computing full matrix vector products. Such high-dimensional FFT-like algorithms, however, typically need to exploit some structure of the sampling pattern and are hence restricted to specific sampling sets.
Even for high-dimensional trigonometric polynomials on a regular full grid, this structure will likely get lost in the dimension-incremental procedure: As soon as previous iterations have identified an unstructured candidate set, this set will no longer have the advantageous structure of the underlying grid that can be exploited for computing the FFT. 

Algorithm~\ref{alg:compute_p_from_Gamma} in contrast, yields a fast transform for arbitrary candidate sets. Consequently, it can be efficiently applied in every iteration no matter how the candidate set looks like, which in turn gives rise to the superior computational complexity of the high dimensional sparse fast Fourier transform resulting from employing Algorithm~\ref{alg:compute_p_from_Gamma} (together with some subsequent low-cost computations) in the role of Algorithm~$\operatorname{A}$ in the context of Algorithm~\ref{algo:sfft_general}.

\subsection{Sample complexity and computational complexity}
\label{sec:complexity}

\subsubsection{Algorithm~\ref{algo:sfft_general} using a general efficient identification Algorithm~$\operatorname{A}$}
\label{sec:complexity:Alg3_AlgA}

Since we proposed to consider Algorithm~\ref{algo:sfft_general} as a meta-algorithm, we start by analyzing its sample complexity as well as its computational complexity in general terms.

First we consider the number of samples used by Algorithm~\ref{algo:sfft_general}
for each step separately. The general results are collected in Table~\ref{tab:sfft_general_complexities}.
In Step~1, we use $r\sum_{t=1}^d K_t\le d\,r\,N_\Gamma$ many sampling values.
For Step~2 we construct $r$ different sampling sets $\mathcal{X}_{t,i}$ for each $t=2,\ldots,d-1$ and, in addition, one sampling set $\mathcal{X}_{d,1}$, i.e., $r\,(d-2)+1$ many.
Since different choices for Algorithm $\operatorname{A}$ require different sampling sets and hence the sampling strategy must be chosen in conjunction with the algorithm substituted for this black box, we denote the sampling sets by $\mathcal{X}_{t,i}=\mathcal{X}_{t,i}(\operatorname{A})$. 
Certainly, the sampling sets $\mathcal{X}_{t,i}(\operatorname{A})$ may depend on the admissible failure probability $\gamma_{\operatorname{A}}$ of Algorithm $\operatorname{A}$.
With this notation, the number of  sampling values used in Step~2 is bounded by
$|\mathcal{X}_{d,1}(\operatorname{A})|+\sum_{t=2}^{d-1}\sum_{i=1}^{r}|\mathcal{X}_{t,i}(\operatorname{A})|$. For Step~3, we need a number of sampling nodes on the order of $\max(s,N_\Gamma)\log (s/\gamma)$ for a parameter $\gamma$. Namely, this is the number of sampling nodes required by \cite[Algorithm~1]{KaPoVo17} to realize a spatial discretization of trigonometric polynomials with frequencies supported in $I^{(1,\ldots,d)}$ with probability at least $1-\gamma$ (see \cite{KaPoVo17} for details). 
\begin{table}[tb]
\centering
\begin{tabular}{lll}
\toprule
& sample complexity & computational complexity\\
\midrule
Step 1 &
$d\,r\,N_\Gamma$&
$d\,r\,N_\Gamma\log{N_\Gamma}$
\\
Step 2&
$|\mathcal{X}_{d,1}(\operatorname{A})|+\sum_{t=2}^{d-1}\sum_{i=1}^{r}|\mathcal{X}_{t,i}(\operatorname{A})|$
&
$d\,r\,\mathcal{C}(\operatorname{A})+d\,r\,s\,\log(rs)$
\\
Step 3 &
$\max(s,N_\Gamma)\log(s/\gamma)$&
$\max(s,N_\Gamma)\log(s/\gamma)(d+\log(s\,N_\Gamma))$
\\
\bottomrule
\end{tabular}
\caption{Sample complexities and computational complexities for the different steps of Algorithm~\ref{algo:sfft_general}, where $\mathcal{C}(\operatorname{A})$ denotes the
maximum computational cost of a single invocation of the efficient identification Algorithm~$\operatorname{A}$ and where the multiple rank-1 lattice approach from \cite[Algorithms~1 and~2]{KaPoVo17} is used in Step~3.}\label{tab:sfft_general_complexities}
\end{table}

Second, we consider the computational complexity, again for each step. In Step~1, we apply $r$ one-dimensional FFTs of maximal size $N_\Gamma$, $d$ times each.
Moreover, each of the $d$ different frequency sets $I^{(t)}$, $t=1,\ldots,d$, is constructed incrementally  in $r$ substeps,
e.g., by sorting vectors $(|\tilde{\hat{p}}_{t,j}|)_{j\in\mathcal{P}_t(\Gamma)}$ of length at most $N_\Gamma$ and handling a sorted vector of length at most $N_\Gamma$.
Accordingly, the computational complexity of Step~1 is in $d\,rN_\Gamma\log N_\Gamma$.

The computational complexity of Step~2 intrinsically depends on the choice of Algorithm~$\operatorname{A}$. This dependency is two-fold. On the one hand, in Step~2c, the algorithm is executed, hence the runtime of this step directly corresponds to the computational complexity of Algorithm~$\operatorname{A}$. On the other hand, different choices for Algorithm~$\operatorname{A}$ will require different sampling sets, which can be costly to construct if specific properties are required. Hence the computational complexity of Step 2a, in which the sampling set is constructed will also depend on Algorithm~$\operatorname{A}$. 
We will subsume the $\operatorname{A}$-dependent contribution to the computational complexity arising in these two steps in a constant $C(\operatorname{A})$, chosen to be a universal upper bound for these  contributions over all possible choices of $t$ and $i$. 
Step 2b also has a mild implicit dependency on Algorithm~$\operatorname{A}$, as it scales with the size of the sampling set. However, it should not be seen as a part of the algorithm, as it also depends on the sampling procedure in the underlying application, which is independent of the algorithm design. In this paper, we will follow the assumption made in most other works on the topic that the sampling procedure is dominated by the other steps of the algorithms, i.e., the computational complexity of Step~2b is also of order $\OO{C(\operatorname{A})}$. At this point, we would remind the reader once again that we tolerate a specific failure probability $\gamma_{\operatorname{A}}$ of Algorithm~$\operatorname{A}$ which may imply a $\gamma_{\operatorname{A}}$ dependence of $C(\operatorname{A})$.

Step~2d is about adding at most $s_\mathrm{local}$ elements to the set $I^{(1,\dots,t)}$; the largest run-time contribution of this step is to avoid listing an existing element again. In   analogy to the considerations of Step~1, we observe a computational complexity in $\OO{d\,r\,s\log(rs)}$ for Step~2d under the assumption $s_\mathrm{local}\sim s$.%

As Step~3 is based on \cite[Algorithms~1 and~2]{KaPoVo17}, the results in this paper provide bounds for the computational complexity of constructing the sampling set and computing the Fourier coefficients. We obtain $\OO{\max(s,N_\Gamma)\log(s/\gamma)(d+\log(s\,N_\Gamma))}$,
with the parameter $\gamma$ as introduced in the beginning of this subsection. 

\subsubsection{Algorithm~\ref{algo:sfft_general} using a modification of  Algorithm~\ref{alg:compute_p_from_Gamma} in the role of Algorithm~$\operatorname{A}$}

\begin{algorithm}[tb]
	\caption{Modification of Algorithm~\ref{alg:compute_p_from_Gamma} limiting the output.}\label{alg:compute_p_from_Gamma_s}
	\begin{tabular}{p{2cm}p{4.0cm}p{8.25cm}}
		Input: 	& $\Gamma\subset\Z^d$ 								& set of frequency candidates\\
		& $\Lambda:=\bigcup_{\ell=1}^L\Lambda(\boldz_\ell,M_\ell)$ 	& sampling nodes, $L$ odd\\
		& $\left\{p(\boldx)\colon\boldx\in\Lambda\right\}$			& sampling values of $p\in\Pi_I$ \\
		& $\tilde{s}$												& sparsity parameter \\
		& $\theta$													& Fourier coefficient threshold \\
	\end{tabular}
	\begin{algorithmic}[1]
		\State Compute $\tilde{I}$ and $(\check{p}_\boldk)_{\boldk\in\tilde{I}}$ using Algorithm~\ref{alg:compute_p_from_Gamma}
		\State Set $\tilde{I}:=\{\boldk\in\tilde{I}\colon: |\check{p}_{\boldk}|\ge\theta\}$ and $\tilde{s}=\min(\tilde{s},|\tilde{I}|)$.
		\State Sort the entries of $(\check{p}_\boldk)_{\boldk\in\tilde{I}}$ by their modulus in descending order to obtain $(\check{p}_{\boldk_1},\check{p}_{\boldk_2}, \dots)$\newline and set $\tilde{I}:=\{\boldk_1,\ldots,\boldk_{\tilde{s}}\}$.
	\end{algorithmic}
	\begin{tabular}{p{2cm}p{4.0cm}p{8.25cm}}
		Output: & $\tilde{I}$						& detected frequency set $|\tilde{I}|\le\tilde{s}$\\
		& $(\check{p}_\boldk)_{\boldk\in\tilde{I}}$	& Fourier coefficients belonging to $\tilde{I}$, i.e., $|\check{p}_\boldk|\ge\theta$%
		\\    \cmidrule{1-3}
		\multicolumn{3}{p{15.1cm}}{Computational complexity: $\OO{L(M\log M+d\,|\Gamma|) + |\Gamma|\log|\Gamma|}$, where $M=\max\{M_\ell\colon \ell=1,\ldots,L\}$}
	\end{tabular}
\end{algorithm}

In this section, we analyze Algorithm~\ref{algo:sfft_general} when using
multiple random rank-1 lattices for constructing the first $t$ components of the sampling sets $\mathcal{X}_{t,i}$ in Step~2a according to Corollary~\ref{cor:prob_bound_reco_trig_poly} and applying Algorithm~\ref{alg:compute_p_from_Gamma} as Algorithm~$\operatorname{A}$ in Step~2c, where we use $J_t$ as the set of frequency candidates~$\Gamma$ for Algorithm~\ref{alg:compute_p_from_Gamma}.
Since we assumed Algorithm~$\operatorname{A}$ as an estimator of the $\tilde{s}$ most significant frequencies of the input signal, we need to apply a slight modification of Algorithm~\ref{alg:compute_p_from_Gamma}
which we summarize in Algorithm~\ref{alg:compute_p_from_Gamma_s}. It is just the application of Algorithm~\ref{alg:compute_p_from_Gamma} with a subsequent additional restriction of the output $\tilde{I}$ guaranteeing that 
this frequency set then fits to the requirements of Step~2c in Algorithm~\ref{algo:sfft_general}.

As above, we assume $s_\mathrm{local}\sim s$, e.g., $s_\mathrm{local}:=2s$, in order to avoid at least this parameter. It suffices to discuss Step~2 since the complexity of the other steps is independent of Algorithm~$\operatorname A$ and has already been investigated in the last section. Some parameters related to the probability of failure and the number of iterations will not be specified in this section, the next section discusses suitable choices. In particular, the failure probability $\gamma_{\operatorname{A}}$ of Algorithm~$\operatorname{A}$ still remains unspecified in this section.
In Table~\ref{tab:sfft_alg2_complexities}, we give an overview of the sample complexities as well as computational complexities of the three steps involving these parameters $r$, $\gamma_{\operatorname{A}}$, and $\gamma$. The corresponding complexities for the suitable parameter choices will again be postponed to the next section.

\begin{table}[tb]
\centering
\begin{tabular}{lll}
\toprule
& sample complexity & computational complexity\\
\midrule
Step 1 &
$d\,r\,N_\Gamma$&
$d\,r\,N_\Gamma\log{N_\Gamma}$
\\
Step 2 (w.h.p.)&
$d\,r\,\max(s,N_\Gamma)\log\frac{s\,N_\Gamma}{\gamma_{\operatorname{A}}}$
&
$d^2\,r\,s\,N_\Gamma\log^2\left(\frac{s\,N_\Gamma}{\gamma_{\operatorname{A}}}\right)$
\\
Step 2 (w.c.)&
$d\,r\,\max(s,N_\Gamma)\log\frac{r\,s\,N_\Gamma}{\gamma_{\operatorname{A}}}$
&
$d^2\,r^2\,s\,N_\Gamma\log^2\left(\frac{r\,s\,N_\Gamma}{\gamma_{\operatorname{A}}}\right)$
\\
Step 3 &
$\max(s,N_\Gamma)\log(s/\gamma)$&
$\max(s,N_\Gamma)\log(s/\gamma)(d+\log(s,N_\Gamma))$
\\
\bottomrule
\end{tabular}
\caption{Sample complexities and computational complexities for the different steps of Algorithm~\ref{algo:sfft_general}, where the efficient identification Algorithm~\ref{alg:compute_p_from_Gamma_s} is used in Step 2
and the multiple rank-1 lattice approach from \cite[Algorithm~1]{KaPoVo17} in Step~3.}
\label{tab:sfft_alg2_complexities}
\end{table}

Again, we start with the sample complexity.
Each of the sampling sets $\mathcal{X}_{t,i}$ is the union of  $\OO{\log(|J_t|/\gamma_{\operatorname{A}})}$ rank-1 lattices in $t$ spatial dimensions, consisting of~$M_{t,i,\ell}=\OO{\max(s,N_\Gamma)}$ lattice nodes each, which are embedded into $d$ dimensions by concatenating all the lattice nodes by $d-t$ fixed components (the same for all nodes in the sampling set $\mathcal{X}_{t,i}$), which are drawn uniformly at random from $\T$.
The generating vectors $\boldz_{t,i,\ell}$  of the $t$-dimensional rank-$1$ lattices are drawn independently of these components and of each other, uniformly at random from $[0,M_{t,i,\ell}-1]^t\cap\Z^t$.
As $|J_t|\lesssim r\,s\,N_\Gamma$, an upper bound of the size of each sampling set is given by $|\mathcal{X}_{t,i}|\lesssim \max(s,N_\Gamma)\,\log(r\,s\,N_\Gamma/\gamma_{\operatorname{A}})$. As the number of these sampling sets  $\mathcal{X}_{t,i}$ is $1+(d-2)\,r$, we obtain that the sampling complexity of Step~2 is $\OO{d\,r\,\max(s,N_\Gamma)\,\log(r\,s\,N_\Gamma/\gamma_{\operatorname{A}})}$.

In order to estimate the computational complexity of Step~2 with Algorithm~\ref{alg:compute_p_from_Gamma_s}
taking the role of Algorithm~$\operatorname{A}$, we distinguish two different cases.
First, we consider the worst case scenario.
For this, we just use the estimate $|J_t|\lesssim r\,s\,N_\Gamma$ and apply Theorem~\ref{thm:intro} for each $t$ and $i$, where  $J_t$ takes the role of the set of frequency candidates~$\Gamma$ in Theorem~\ref{thm:intro}. Accordingly, we observe a computational complexity in 
\begin{equation}
\OO{d\,r\,\Big(\max(s,N_\Gamma)\log(s\,N_\Gamma)+d\,r\,s\,N_\Gamma\Big)\log\frac{r\,s\,N_\Gamma}{\gamma_{\operatorname{A}}}}.
\label{eq:wc_complexity_alg3_alg1}
\end{equation}
This estimate, however, is far from order optimal in many cases, as the bound on $|J_t|$ is approximately tight only when the sets $\tilde{J}_{t-1,i}$ constructed for different values of $i$ have very little overlap. 

In the event that
Step~1 of Algorithm~\ref{algo:sfft_general}, i.e., $P_t(I)\subset I^{(t)}$ and, in addition, all instances of Algorithm~\ref{alg:compute_p_from_Gamma_s} are successful (in the sense of Corollary~\ref{cor:prob_bound_reco_trig_poly}, cf.\ Section~\ref{ssec:parameter_alg4_as_algA_in_alg3}), these sets will have large overlap.
As we will discuss in the next section, this event has high probability for suitable parameter choices. That is, in such scenarios, the computational complexity will be considerably smaller with high probability.
More precisely, in that event one has $|J_t|\le sN_\Gamma$,  so one obtains a computational complexity that is reduced by a factor of $r$. This yields a computational complexity in
\begin{equation}
\OO{d\,r\,(\max(s,N_\Gamma)\log(s\,N_\Gamma)+d\,s\,N_\Gamma)\log(s\,N_\Gamma/\gamma_{\operatorname{A}})}.
\label{eq:whp_complexity_alg3_alg1}
\end{equation}
A similar argument shows that with high probability, one also observes a slight improvement of the sample complexity in the sense that the 
logarithmic term will no longer depend on~$r$.

\subsubsection{Choosing the parameters in Algorithm~\ref{algo:sfft_general} for Algorithm~\ref{alg:compute_p_from_Gamma_s} in the role of Algorithm~$\operatorname{A}$}
\label{ssec:parameter_alg4_as_algA_in_alg3}

Up to now, we have not discussed how to choose the parameters $r$, $\gamma_{\operatorname{A}}$, and $\gamma$. Following the probability estimates from \cite[Lemma 4.4]{KaPoVo17} motivates the choice $r:=\lceil 2\,s\,\log\left(\frac{3\,d\,s}{\delta}\right)\rceil$.
This choice ensures that each of the non-zero Fourier coefficients $\hat{p}_\boldk$ of $p$ as in \eqref{eq:trig_poly_p} can be  detected in Step~1 as well as Step~2 with high probability.
The key idea of \cite{KaPoVo17} that we are also exploiting here, is that fixing the last $d-t$ components of all vectors in the sampling set to the same fixed values, allows for the estimation of a certain projection of the Fourier coefficients. Repeating this process $r$ times for different choices of the $d-t$ components ensures that with a probability of at least $1-\tfrac{\delta}{3\,d}$, each of the active Fourier coefficients
is projected to some coefficient that is not less than $\theta$ at least once among the different projections.
Since we take the union of the frequency sets $\tilde{J}_{t,i}$, $i=1,\ldots,\tilde{r}$, in Step~2d,
it is sufficient to only regard the frequencies that belong to coefficients of at least $\theta$ in modulus
in Step~2c.
Clearly, Algorithm~\ref{alg:compute_p_from_Gamma} will compute all projections of the active Fourier coefficients, and in particular those that are at least $\theta$, with failure probabilities estimated in Corollary~\ref{cor:prob_bound_reco_trig_poly}. Accordingly, the failure probabilities for detecting all the frequencies with coefficients not less than $\theta$ are also bounded by these estimates and
as a consequence, we can apply the estimates on the failure probabilities in Corollary~\ref{cor:prob_bound_reco_trig_poly} to Algorithm~\ref{alg:compute_p_from_Gamma_s} since $\tilde{s}\ge s$ is fulfilled.

That is why we set $\gamma_{\operatorname{A}}:=\delta/(3\,d\,r)$ in Step~2, which entails that for each fixed $t$ and fixed $i$, Algorithm~\ref{alg:compute_p_from_Gamma_s} has a failure probability of at most $\delta/(3\,d\,r)$
for detecting all frequencies belonging to Fourier coefficients of at least $\theta$.
Similarly, we fix $\gamma:=\delta/(3\,d)$ for Step~3.

In analogy to \cite[Theorem~4.6]{KaPoVo17}, the total failure probability can now be estimated via a union bound over the different parts.
We obtain failure probabilities of $\delta/3$ for Step~1, $(d-2)\delta/(3\,d)+((d-2)r+1)\delta/(3\,d\,r)$ for Step~2, and $\delta/(3\,d)$ for Step~3.
Accordingly, the total failure probability is less than $\delta$. 

\begin{table}[tb]
\centering
\begin{tabular}{lll}
\toprule
& sample complexity & computational complexity\\
\midrule
Step 1 &
$d\,s\,N_\Gamma\,\log{\frac{d\,s}{\delta}}$&
$d\,s\,N_\Gamma\,\log^2{\frac{d\,s\,N_\Gamma}{\delta}}$
\\
\rule{0em}{1.4em}Step 2 (w.h.p.)&
$d\,s\max(s,N_\Gamma)\log^2\frac{d\,s\,N_\Gamma}{\delta}$
&
$d^2s^2N_\Gamma\log^3\frac{d\,s\,N_\Gamma}{\delta}$
\\
\rule{0em}{1.4em}Step 2 (w.c.)&
$d\,s\max(s,N_\Gamma)\log^2\frac{d\,s\,N_\Gamma}{\delta}$
&
$d^2s^3N_\Gamma\log^3\frac{d\,s\,N_\Gamma}{\delta}$
\\
\rule{0em}{1.4em}Step 3 &
$\max(s,N_\Gamma)\log(d\,s/\delta)$&
$\max(s,N_\Gamma)\log(d\,s/\delta)(d+\log(s,N_\Gamma))$
\\
\bottomrule
\end{tabular}
\caption{Sample complexities and computational complexities for the different steps of Algorithm~\ref{algo:sfft_general}, where the efficient identification Algorithm~\ref{alg:compute_p_from_Gamma_s} is used in Step 2
and $r:=\ceil{2\,s\,\ln{\frac{3\,d\,s}{\delta}}}$ is chosen in line with \cite[Theorem~4.6]{KaPoVo17}.
Step 3 is realized via the multiple rank-1 lattice approach of \cite[Algorithm~1]{KaPoVo17}.}
\label{tab:sfft_whp_complexities}
\end{table}

\begin{sloppypar}
The resulting sample complexities and computational complexities for this choice of parameters are summarized in Table~\ref{tab:sfft_whp_complexities}. For Step 2 we list a worst case (w.c.) bound as well as a bound that hold with a high probability of at least $1-\delta$ (w.h.p.).

\end{sloppypar}

For the computational complexity of Step~2,
we obtain simpler expressions by bounding the worst case estimate via
\begin{equation*}
d\,s\,\left(\log\frac{d\,s}{\delta}\right)\,\left(\max(s,N_\Gamma)\log(s\,N_\Gamma)+d\,s^2\,N_\Gamma\,\log\frac{d\,s}{\delta}\right)\log\frac{d\,s\,N_\Gamma}{\delta}\lesssim d^2s^3N_\Gamma\log^3\frac{d\,s\,N_\Gamma}{\delta}
\end{equation*}
and the high probability estimate via
\begin{equation*}
d\,s\,\left(\log\frac{d\,s}{\delta}\right)\,(\max(s,N_\Gamma)\log(s\,N_\Gamma)+d\,s\,N_\Gamma)\log\frac{d\,s\,N_\Gamma}{\delta}
\lesssim d^2 s^2 N_\Gamma\log^3\frac{d\,s\,N_\Gamma}{\delta}.
\end{equation*}

\section{Numerical results}

In this section, we validate our theoretical findings by numerical experiments. We first demonstrate the feasibility for unstructured candidate sets by choosing both the candidate set and the set of active frequencies at random. Second, as an example of highly structured candidate sets, we investigate the hyperbolic cross. Lastly we study the performance of our method in the context of a dimension-incremental sparse FFT in high dimensions.

\subsection{Sparse FFT for arbitrary candidate sets as introduced in Section~\ref{sec:main_result}}
\label{ssec:num1}

We start by validating our approach in the framework of Section~\ref{sec:main_result}.
That is, we consider different frequency candidate sets $\Gamma\subset\Z^d$, $|\Gamma|<\infty$,  and apply our method to recover trigonometric polynomials $p$, cf.~\eqref{eq:trig_pol}, with frequencies supported on small subsets $I\subset\Gamma$. 

\subsubsection{Identifying potential false positives and potential false negatives}\label{sec:PFPN}

As explained above, a main advantage of our recovery guarantees as compared to other results with similar sample complexity is that we do not assume a random signal model, but we show recovery with high probability for an arbitrary signal. Even stronger, we only analyze aliasing properties of the support, cf.\ formula~\eqref{eq:r1l_aliasing_formula} and the related discussion in Section~\ref{sec:main_result:background}. So for fixed support $I$ we guarantee that with high probability all signals supported on $I$ can be recovered. 

In our numerical simulations, we also aim to illustrate this strong property.
For that, we employ a worst case measure for the support detection. 
For a given support $I$, we compute all the {\em potential} false positives (PFP) and {\em potential} false negatives (PFN), that is, all index vectors $\mathbf j$ that arise as an alias of some $\boldk\in I$. When $\boldj\in I$, a cancellation of the true coefficient and the aliased coefficient can have the effect that the associated frequency is not detected as active. For $\boldj\notin I$, the aliased coefficient will result in an unjustified detection of the coefficient. Especially the potential false negatives will only be realized for very specific choices of coefficients (that would be unlikely under a random model), and it is not clear if there actually are coefficients that realize multiple of them simultaneously. For this reason, the bounds on the success rates we empirically compute are somewhat conservative, but they certainly form a lower bound for the true rates. %

More precisely, our empirical evaluation is based on the observation
that to decide if a given frequency $\boldk$ is in $I$, our method
considers a set of measurements, each of which constitutes the sum over
the true (potentially zero if $\boldk \notin I$) coefficient
corresponding to $\boldk$ and all the aliased coefficients with respect
to some random rank-1 lattice, and classifies $\boldk$ as a member of
$I$ if fewer than half of these measurements are zero. Thus a sufficient
condition to prevent that $\boldk$ is wrongly classified as a member or
not a member of $I$ is that for more than half of these random lattices,
$\boldk$ aliases to no other $\boldj \in I$.

To count how many $\boldj \in I$ the frequency $\boldk$ aliases to, we use the following trick: We apply the method to the auxiliary polynomial 
\begin{equation}
p(\boldx)=\sum_{\boldk\in I}\e^{2\pi\ii\boldk\cdot\boldx}.
\label{eq:eq:dirichlet_I}
\end{equation}
As all the coefficients indexed by $I$ are $1$, the aforementioned sum over the aliased coefficients for a realization of the random lattice is nothing but a counting measure applied to the intersection of $I$ and the set of indices that $\boldk$ aliases to. Consequently, when no aliasing happens, the sum will be $0$ for $\boldk \notin I$ and $1$ for $\boldk\in I$. Both of these values are the smallest possible, so if the median $\check p_\boldk$ in~\eqref{eq:def_median_check_p} takes this value, aliasing happens in fewer than half of the cases, as desired, and hence both false positives and false negatives are excluded for any choice of values of the active coefficients. Otherwise, we count $\boldk$ as a potential false positive (PFP) or potential false negative (PFN), respectively. This method to empirically identify the potential false detections is illustrated in Table~\ref{tab:wc_scenarios}.

\begin{table}[H]
\begin{center}
\begin{tabular}{|l|l|c|c|c|}
	\hline
	\multicolumn{2}{|c|}{}& 	\multicolumn{3}{|c|}{Coefficient recovered by Algorithm \ref{alg:compute_p_from_Gamma}} \\
\cline{3-5}
 \multicolumn{2}{|c|}{}& $\check{p}_\boldk=0$ & $\check{p}_\boldk=1$& $\check{p}_\boldk\geq 2$\\
\hline
\multirow{2}{2.4cm}{{\bfseries true} Fourier coefficient of $p$} & $\hat{p}_\boldk=1$ (i.e., $\boldk\in I$) & not possible & not a PFN & PFN\\
\cline{2-5}
&$\hat{p}_\boldk=0$ (i.e., $\boldk\in \Gamma\setminus I$)  & not a PFP  & \multicolumn{2}{c|}{PFP}\\
\hline
\end{tabular}
\end{center}
\caption{Empirical detection of potential false negatives (PFN) and potential false positives (PFP) for known $\Gamma$ and $I$ via Algorithm~\ref{alg:compute_p_from_Gamma} applied to the
auxiliary trigonometric polynomial given in \eqref{eq:eq:dirichlet_I}}\label{tab:wc_scenarios}
\end{table}

\begin{Remark}\label{rem:do_not_use_r1l_info_for_pfn}
On the one hand, the empirical detection of the PFNs and PFPs via Algorithm~\ref{alg:compute_p_from_Gamma} applied to the
auxiliary trigonometric polynomial~\eqref{eq:eq:dirichlet_I} will always give correct results. On the other hand, when additionally using the postprocessing discussed in Remark~\ref{rem:use_r1l_info_or_ls}, we may obtain ``wrong'' counts if there are any potential false negatives present. The reason for this is that the PFN frequencies will be in the output~$\tilde{I}$ of Algorithm~\ref{alg:compute_p_from_Gamma} and consequently, taken into consideration by the postprocessing discussed in Remark~\ref{rem:use_r1l_info_or_ls}. Whenever the set of used rank-1 lattices provides a spatial discretization of the space $\Pi_{\tilde{I}}$ of trigonometric polynomials, all PFNs will be filtered out. In reality, however, it may happen that a potential false negative does not appear in the output~$\tilde{I}$ at all if aliasing Fourier coefficients cancel out each other, and consequently, the postprocessing could possibly yield incorrect Fourier coefficients and filter out energetic frequencies. As the cancellation strongly depends on the used rank-1 lattices and the Fourier coefficients of the function under consideration, the cancellation should be unlikely in practice.
\end{Remark}

\subsubsection{Random frequency sets and candidate sets in three dimensions %
	}\label{sec:num_sec2_random}

Here, we investigate the accuracy of multiple random rank-1 lattice sampling for the exactly sparse case.
We consider the reconstruction of the Fourier coefficients~$\hat{p}_\boldk$ of three-variate trigonometric polynomials~$p$ of sparsity $|I|=1\,000$.
For each $L\in\{9,11,13,\ldots,37\}$, we fix a multiple random rank-1 lattice configuration $\Lambda(\boldz_1,M_1,\ldots,\boldz_L,M_L)$, where each single \hbox{rank-1} lattice is of size $M_\ell=10\,331>10.33\cdot 1\,000$.
Now we repeat  the following computations 1\,000 times.
We choose an index set of possible frequencies $\Gamma\subset\{-1\,000,-999,\ldots,1\,000\}^3$, $|\Gamma|=10^7$, and the index set of active frequencies $I\subset\Gamma$, $|I|=1\,000$, both uniformly at random. For these choices of $\Gamma$ and $I$ we identify potential false positives (PFP) and potential false negatives (PFN) as explained in Section~\ref{sec:PFPN} and illustrated in Table~\ref{tab:wc_scenarios}. 

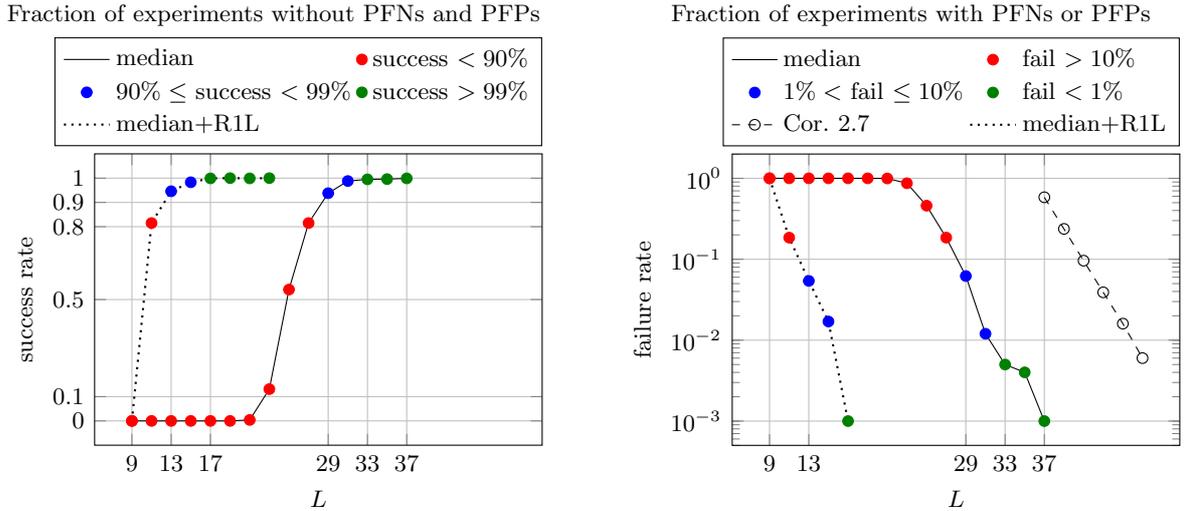
\begin{figure}[bth]
		\begin{tikzpicture}[baseline]
		\begin{axis}[
			font=\footnotesize,
			enlarge x limits=true,
			enlarge y limits=true,
			height=0.35\textwidth,
			grid=major,
			width=0.48\textwidth,
			xtick={9,13,17,29,33,37},
			ytick={0,0.1,0.5,0.8,0.9,1},
            xmax=47,
			ymin=0,ymax=1,
			xlabel={$L$},
			ylabel={success rate},
			legend style={legend cell align=left, at={(1,1.4)}}, %
			legend columns = 2,
			title={Fraction of experiments without PFNs and PFPs},
            title style={at={(0.4,1.35)}},
		]
		\addplot[no marks,black] coordinates {
 (9,0.000) (11,0.000) (13,0.000) (15,0.000) (17,0.000) (19,0.000) (21,0.004) (23,0.131) (25,0.541) (27,0.815) (29,0.938) (31,0.988) (33,0.995) (35,0.996) (37,0.999)
		};
		\addlegendentry{median}
		\addplot[only marks,red,mark=*,mark options={solid}] coordinates {
 (9,0.000) (11,0.000) (13,0.000) (15,0.000) (17,0.000) (19,0.000) (21,0.004) (23,0.131) (25,0.541) (27,0.815)
		};
		\addlegendentry{success $< 90\%$}
		\addplot[only marks,blue,mark=*,mark options={solid}] coordinates {
 (29,0.938) (31,0.988)
		};
		\addlegendentry{$90\%\leq$ success $<99\%$}
		\addplot[darkgreen,only marks,mark=*,mark options={solid}] coordinates {
 (33,0.995) (35,0.996) (37,0.999)
		};
		\addlegendentry{success $> 99\%$}
		\addplot[dotted,thick,no marks,black] coordinates {
 (9,0.000) (11,0.815) (13,0.946) (15,0.983) (17,0.999) (19,1.000) (21,0.999) (23,1.000)
		};
		\addlegendentry{median+R1L}
		\addplot[only marks,forget plot,red,mark=*,mark options={solid}] coordinates {
 (9,0.000) (11,0.815)
		};
		\addplot[only marks,forget plot,blue,mark=*,mark options={solid}] coordinates {
 (13,0.946) (15,0.983)
		};
		\addplot[only marks,forget plot,darkgreen,mark=*,mark options={solid}] coordinates {
 (17,0.999) (19,1.000) (21,0.999) (23,1.000)
		};
		\end{axis}
		\end{tikzpicture}
\hfill
		\begin{tikzpicture}[baseline]
		\begin{axis}[
			font=\footnotesize,
			enlarge x limits=true,
			enlarge y limits=true,
			height=0.35\textwidth,
			grid=major,
			width=0.48\textwidth,
			xtick={9,13,29,33,37},
			xlabel={$L$},
			ylabel={failure rate},
			legend style={legend cell align=left, at={(1,1.4)}}, %
			legend columns = 2,
			ymode=log,
			title={Fraction of experiments with PFNs or PFPs},
            title style={at={(0.4,1.35)}},
		]
		\addplot[no marks,black] coordinates {
 (9,1.000) (11,1.000) (13,1.000) (15,1.000) (17,1.000) (19,1.000) (21,0.996) (23,0.869) (25,0.459) (27,0.185) (29,0.062) (31,0.012) (33,0.005) (35,0.004) (37,0.001)
		};
		\addlegendentry{median}
		\addplot[red,only marks,mark=*,mark options={solid}] coordinates {
 (9,1.000) (11,1.000) (13,1.000) (15,1.000) (17,1.000) (19,1.000) (21,0.996) (23,0.869) (25,0.459) (27,0.185)
		};
		\addlegendentry{fail $>10\%$}
		\addplot[blue,only marks,mark=*,mark options={solid}] coordinates {
 (29,0.062) (31,0.012)
		};
		\addlegendentry{$1\%<$ fail $\leq 10\%$}
		\addplot[darkgreen,only marks,mark=*,mark options={solid}] coordinates {
 (33,0.005) (35,0.004) (37,0.001)
		};
		\addlegendentry{fail $< 1\%$}
		\addplot[black,dashed,mark=o,mark options={solid}] coordinates {
 (37,0.583) (39,0.237) (41,0.096) (43,0.039) (45,0.016) (47,0.006)
		};
		\addlegendentry{Cor.~\ref{cor:prob_bound_reco_trig_poly}}
		\addplot[no marks,black,dotted,thick] coordinates {
 (9,1.000) (11,0.185) (13,0.054) (15,0.017) (17,0.001) (19,0.000) %
		};
		\addlegendentry{median+R1L}
		\addplot[only marks,forget plot,red,mark=*,mark options={solid}] coordinates {
 (9,1.000) (11,0.185)
        };
		\addplot[only marks,forget plot,blue,mark=*,mark options={solid}] coordinates {
 (13,0.054) (15,0.017)
		};
		\addplot[only marks,forget plot,darkgreen,mark=*,mark options={solid}] coordinates {
 (17,0.001) (19,0.000) %
		};
		\end{axis}
		\end{tikzpicture}
		\caption{Success rate and failure rate with respect to the number~$L$ of used rank\mbox{-}1 lattices for random frequency sets $I$ and $\Gamma$, $|I|=1000$ and $|\Gamma|=10^7$.
        Solid lines correspond to Algorithm~\ref{alg:compute_p_from_Gamma}, the dotted lines corresponds to Algorithm~\ref{alg:compute_p_from_Gamma} with postprocessing as discussed in Remark~\ref{rem:use_r1l_info_or_ls}, and the dashed line shows the theoretical bounds from Corollary~\ref{cor:prob_bound_reco_trig_poly}.}
        \label{fig:success_rate_L_random_gamma_10m}
		\end{figure}

In Figure~\ref{fig:success_rate_L_random_gamma_10m}, we plot the success and failure rates for these experiments as solid lines and filled circles for various choices for the number $L$ of rank-1 lattices used in the configuration. Here we count an experiment as a success when no PFPs and PFNs are identified.  

 To put these success rate into perspective, recall that by 
 Corollary~\ref{cor:prob_bound_reco_trig_poly}, the probability that Algorithm~\ref{alg:compute_p_from_Gamma} will correctly determine a trigonometric polynomial using the multiple random rank-1 lattice construction is bounded from below by $\max(0,1-10^7\cdot 9.331^{-8.331/41.324\, L})$. This number is positive for odd $L\ge 37$;  %
 for $L=41$ one already obtains a lower bound on the success probability of $0.90$.
As expected some estimates in the proof of the corollary are not tight. On the one hand, for a similar target failure rate, the number $L$ of used rank-1 lattices in the numerical tests is lower by approximately one third compared to the corresponding theoretical bound in Corollary~\ref{cor:prob_bound_reco_trig_poly}. On the other hand, for $L=37$ rank-1 lattices, the empirical failure rate behaves distinctly better with a value of only 0.001 compared to the higher theoretical bound by Corollary~\ref{cor:prob_bound_reco_trig_poly} of 0.583.
To illustrate this difference, we include our bound for the failure rate given by Corollary~\ref{cor:prob_bound_reco_trig_poly} for comparison as dashed line and unfilled circles. We observe that despite large constants, the exponential decay of both theory and experiments match well.

Furthermore, we observe that the empirical failure rate is less than 0.01 for $L\ge 33$. For $L=33$, we require 340\,891 samples, which is only $\approx 1/29$ of the samples required for a full discrete Fourier transform on $\Gamma$ and only $\approx 1/23\,500$ of the samples of a fast Fourier transform on $\{-1\,000,-999,\ldots,1\,000\}^3$. Another remarkable observation is that the success probability increases from less than 0.1 to more than 0.9 within a small range of $L$, which also reflects the logarithmic dependence of $L$ on the failure probability.

\begin{figure}[bth]
		\begin{tikzpicture}[baseline]
		\begin{axis}[
			font=\footnotesize,
			enlarge x limits=true,
			enlarge y limits=true,
			height=0.35\textwidth,
			grid=major,
			width=0.48\textwidth,
			xtick={13,15,19,37},
			ytick={0,1,2,3},
			ymin=0,ymax=3,
			xlabel={$L$},
			ylabel={max. \# (in 1\,000 runs)},
			legend style={legend cell align=left}, legend pos=south east,
			legend columns = 1,
			title={Maximal number of PFNs},
		]
		\addplot[only marks,red,mark=square*,mark options={solid}] coordinates {
 (11,3) (13,2) (15,2) (17,1) (21,1) (27,1)
		};
		\addplot[only marks,darkgreen,mark=square*,mark options={solid}] coordinates {
 (19,0) (23,0) (25,0) (29,0) (31,0) (33,0) (35,0) (37,0)
		};
		\end{axis}
		\end{tikzpicture}
\hfill
		\begin{tikzpicture}[baseline]
		\begin{axis}[
			font=\footnotesize,
			enlarge x limits=true,
			enlarge y limits=true,
			height=0.35\textwidth,
			grid=major,
			width=0.48\textwidth,
			xtick={17,19,23,37},
			ytick={1,10,100},
			ymin=0,ymax=100,
			xlabel={$L$},
			ylabel={max. \# (in 1\,000 runs)},
			legend style={legend cell align=left}, legend pos=south east,
			legend columns = 1,
			title={Maximal number of PFPs},
		]
		\addplot[red,no marks,mark options={solid}] coordinates {
 (15,230) (17,88)
		};
		\addplot[black,mark=square*,mark options={solid}] coordinates {
 (17,88) (19,36) (21,16) (23,8)
		};
		\addplot[blue,mark=square*,mark options={solid}] coordinates {
 (23,8) (25,5) (27,3) (29,3) (31,1) (33,1) (35,1) (37,1)
		};
		\end{axis}
		\end{tikzpicture}

		\begin{tikzpicture}[baseline]
		\begin{axis}[
			font=\footnotesize,
			enlarge x limits=true,
			enlarge y limits=true,
			height=0.35\textwidth,
			grid=major,
			width=0.48\textwidth,
			xtick={15,19,21,23,37},
			ytick={0,0.1,0.5,0.8,0.9,1},
			ymin=0,ymax=1,
			xlabel={$L$},
			ylabel={success rate},
			legend style={legend cell align=left}, legend pos=south east,
			legend columns = 1,
			title={Modified success rate allowing for
				$\leq$ $10$ PFPs},
            title style={at={(0.4,1)}},
		]
		\addplot[black,forget plot,no marks] coordinates {
 (15,0.000 )(17,0.000) (19,0.017) (21,0.950) (23,1.000) (25,1.000) (27,0.999) (29,1.000) (31,1.000) (33,1.000) (35,1.000) (37,1.000)
		};
		\addplot[only marks,red,mark=*,mark options={solid}] coordinates {
 (15,0.000 )(17,0.000) (19,0.017)
		};
		\addlegendentry{success $< 90\%$}
		\addplot[only marks,blue,mark=*,mark options={solid}] coordinates {
 (21,0.950)
		};
		\addlegendentry{$90\% \leq$ success $\leq 99\%$}
		\addplot[only marks,darkgreen,mark=*,mark options={solid}] coordinates {
(23,1.000) (25,1.000) (27,0.999) (29,1.000) (31,1.000) (33,1.000) (35,1.000) (37,1.000)
		};
		\addlegendentry{success $> 99\%$}
		\end{axis}
		\end{tikzpicture}
\hfill
		\begin{tikzpicture}[baseline]
		\begin{axis}[
			font=\footnotesize,
			enlarge x limits=true,
			enlarge y limits=true,
			height=0.35\textwidth,
			grid=major,
			width=0.48\textwidth,
			xtick={15,17,37},
			ytick={0,0.1,0.5,0.8,0.9,1},
			ymin=0,ymax=1,
			xlabel={$L$},
			ylabel={success rate},
			legend style={legend cell align=left}, legend pos=south east,
			legend columns = 1,
			title={Modified success rate allowing for
					$\leq$ $100$ PFPs}, %
            title style={at={(0.4,1)}},
		]
		\addplot[forget plot,no marks,black] coordinates {
 (15,0.000) (17,0.999) (19,1.000) (21,0.999) (23,1.000) (25,1.000) (27,0.999) (29,1.000) (31,1.000) (33,1.000) (35,1.000) (37,1.000)
		};
		\addplot[only marks,red,mark=*,mark options={solid}] coordinates {
 (15,0.000)
		};
		\addlegendentry{success $< 90\%$}
		\addplot[only marks,darkgreen,mark=*,mark options={solid}] coordinates {
 (17,0.999) (19,1.000) (21,0.999) (23,1.000) (25,1.000) (27,0.999) (29,1.000) (31,1.000) (33,1.000) (35,1.000) (37,1.000)
		};
		\addlegendentry{success $> 99\%$}
		\end{axis}
		\end{tikzpicture}
\caption{Additional performance measures for the random frequency sets $I$ and $\Gamma$ from Figure~\ref{fig:success_rate_L_random_gamma_10m} for Algorithm~\ref{alg:compute_p_from_Gamma}.\\
	{\em First row:}  Maximal numbers of aliasing frequencies within $I$ (PFNs) and aliasing frequencies in $\Gamma\setminus I$ with frequencies within $I$ (PFPs).\\
	{\em Second row:} Modified success rates with no potential false negatives, but some potential false positives allowed. }\label{fig:success_rate_L_random_detail_gamma_10m}
\end{figure}

In Figure~\ref{fig:success_rate_L_random_detail_gamma_10m}, we consider the aliasing effects in more detail. For $L\ge 29$, we observe no potential false negatives, which means that the output $\tilde{I}$ of Algorithm~\ref{alg:compute_p_from_Gamma} contains all frequencies that belong to non-zero Fourier coefficients, i.e., $I\subset\tilde{I}$.
However, %
 for $20\leq L\leq 37$,
 we still observe that $\tilde{I}\neq I$ in some test runs due to a small number of potential false positives. For our goal of identifying the active frequencies and their associated coefficients, this is much less severe than a false negative would be. 
 
 This observation directly yields a refined definition of success. Namely, we consider Algorithm~\ref{alg:compute_p_from_Gamma} to be successful if no potential false negative and at most a predefined number of potential false positives are observed. The lower plots in 
Figure~\ref{fig:success_rate_L_random_detail_gamma_10m} visualize this refined success rate when this predefined number is chosen to be 10 and 100, respectively.
We observe success in more than $99\%$ of the test runs for $L\ge 23$ and $L\ge 17$ for at most ten potential false positives and for at most 100, respectively.
Furthermore, we note that the success rate seems to exhibit a sharp phase transition, increasing
from less than 0.02 to more than 0.95 from one odd $L$ to the next.

Without providing any further details, we point out that even for significantly smaller $L\ge 17$, we observe at most a small number of false negatives and at most a small number of false positives, which could be a good reason for iterative applications of Algorithm~\ref{alg:compute_p_from_Gamma} with relatively small numbers $L$ of rank-1 lattices and successively reduced frequency support $I$.

As an alternative to this idea, one can additionally use the available rank-1 lattice information to filter the false positives in a postprocessing step as discussed in Remark~\ref{rem:use_r1l_info_or_ls}. In accordance with Remark~\ref{rem:do_not_use_r1l_info_for_pfn}, we only apply the postprocessing if there are no potential false negatives. The corresponding success rates are plotted in Figure~\ref{fig:success_rate_L_random_gamma_10m} as dotted lines and filled circles. We observe a distinct improvement having success rates of more than $99\%$ already for $L\geq 17$.

As briefly indicated in Remark~\ref{rem:use_r1l_info_or_ls}, one could also determine the Fourier coefficients by solving the linear system arising from the restriction to all identified frequencies (including false negatives). Already for minor oversampling, this system will often have a unique solution -- with zero coefficients associated to the false positives. Naturally this approach will require that there are no false negatives and only a small number of false positives.

\subsubsection{Weighted and unweighted hyperbolic crosses as frequency sets and candidate sets in eight dimensions}

In the following, we investigate the detection accuracy for deterministic, structured frequency index sets. To this end, we fix an unweighted eight-dimensional hyperbolic cross
$$\Gamma:=\{\boldk\in \Z^8\colon\prod_{t=1}^8\max(1,|k_t|)\le 32\}\subset[-32,32]^8,$$
$|\Gamma|=10\,665\,297$, as a frequency candidate set. Furthermore, we fix the set of active frequencies $I\subset\Gamma$ of cardinality $|I|=1\,069$ given by
$$
I:=\{\boldk\in\Z^8\colon \prod_{t=1}^8 \max(1,t^{1.08} |k_t|) \leq 32 \}\subset\Gamma,
$$
which is  an eight-dimensional weighted hyperbolic cross.

For each $L\in\{9,11,13,\ldots,37\}$, we repeat the following procedure 1\,000 times: We draw a multiple random rank-1 lattice configuration $\Lambda(\boldz_1,M_1,\ldots,\boldz_L,M_L)$, where each single rank-1 lattice is of size $M_\ell=11\,047>10.33\cdot 1\,069$, and we identify potential false positives (PFP) and potential false negatives (PFN) using Algorithm~\ref{alg:compute_p_from_Gamma} as explained in Section~\ref{sec:PFPN} and illustrated in Table~\ref{tab:wc_scenarios}. %

The plots in Figures~\ref{fig:success_rate_L_hc_weightedhc} and~\ref{fig:success_rate_L_hc_weightedhc_detail} are analogous to those in Section~\ref{sec:num_sec2_random}.
We also observe a similar behavior of the success rates and the failure rates, in line with the decay rates predicted by Corollary~\ref{cor:prob_bound_reco_trig_poly}.
For $L= 31$ -- this corresponds to $342\,427$ samples -- we achieve a success rate of more than 0.99, if we do not allow for any PFNs and PFPs. For the modified notion of success with no PFNs, but 10 or 100 PFPs allowed, a success rate of more than 0.99 can be achieved already for $L= 23$ or $L=19$, respectively.

As before,
if one additionally uses the available rank-1 lattice information, cf.\ Remark~\ref{rem:use_r1l_info_or_ls} and~\ref{rem:do_not_use_r1l_info_for_pfn}, then we observe success rates of more than $99\%$ already for $L\geq 15$, see the dotted lines and filled circles in Figure~\ref{fig:success_rate_L_hc_weightedhc}.

\begin{figure}[tbh]
		\begin{tikzpicture}[baseline]
		\begin{axis}[
			font=\footnotesize,
			enlarge x limits=true,
			enlarge y limits=true,
			height=0.35\textwidth,
			grid=major,
			width=0.48\textwidth,
			xtick={9,13,27,31,37},
			ytick={0,0.1,0.5,0.8,0.9,1},
			ymin=0,ymax=1,
			xlabel={$L$},
			ylabel={success rate},
			legend style={legend cell align=left, at={(1,1.4)}}, %
			legend columns = 2,
			title={Fraction of experiments without PFNs or PFPs},
            title style={at={(0.4,1.35)}},
		]
		\addplot[no marks,black] coordinates {
 (9,0.000) (11,0.000) (13,0.000) (15,0.000) (17,0.000) (19,0.002) (21,0.071) (23,0.442) (25,0.741) (27,0.902) (29,0.965) (31,0.993) (33,0.995) (35,0.997) (37,1.000)
		};
		\addlegendentry{median}
		\addplot[red,only marks,mark=*,mark options={solid}] coordinates {
 (9,0.000) (11,0.000) (13,0.000) (15,0.000) (17,0.000) (19,0.002) (21,0.071) (23,0.442) (25,0.741)
		};
		\addlegendentry{success $< 90\%$}
		\addplot[blue,only marks,mark=*,mark options={solid}] coordinates {
 (27,0.902) (29,0.965)
		};
		\addlegendentry{$90\% \leq$ success $\leq 99\%$}
		\addplot[darkgreen,only marks,mark=*,mark options={solid}] coordinates {
 (31,0.993) (33,0.995) (35,0.997) (37,1.000)
		};
		\addlegendentry{success $> 99\%$}
		\addplot[no marks,black,dotted,thick] coordinates {
 (9,0.000) (11,0.954) (13,0.982) (15,0.995) (17,0.996) (19,1.000) (21,1.000) (23,1.000)
		};
		\addlegendentry{median+R1L}
		\addplot[only marks,forget plot,red,mark=*,mark options={solid}] coordinates {
 (9,0.000)		};
		\addplot[only marks,forget plot,blue,mark=*,mark options={solid}] coordinates {
 (11,0.954) (13,0.982)
		};
		\addplot[only marks,forget plot,darkgreen,mark=*,mark options={solid}] coordinates {
 (15,0.995) (17,0.996) (19,1.000) (21,1.000) (23,1.000)
		};
		\end{axis}
		\end{tikzpicture}
\hfill
		\begin{tikzpicture}[baseline]
		\begin{axis}[
			font=\footnotesize,
			enlarge x limits=true,
			enlarge y limits=true,
			height=0.35\textwidth,
			grid=major,
			width=0.48\textwidth,
			xtick={9,13,27,31,37},
            ymax=1,
			xlabel={$L$},
			ylabel={failure rate},
			legend style={legend cell align=left, at={(1,1.4)}}, %
			legend columns = 2,
			ymode=log,
			title={Fraction of experiments with PFNs or PFPs},
            title style={at={(0.4,1.35)}},
		]
		\addplot[black,no marks] coordinates {
 (9,1.000) (11,1.000) (13,1.000) (15,1.000) (17,1.000) (19,0.998) (21,0.929) (23,0.558) (25,0.259) (27,0.098) (29,0.035) (31,0.007) (33,0.005) (35,0.003) (37,0.000)
		};
        \addlegendentry{median}
		\addplot[red,only marks,mark=*,mark options={solid}] coordinates {
 (9,1.000) (11,1.000) (13,1.000) (15,1.000) (17,1.000) (19,0.998) (21,0.929) (23,0.558) (25,0.259)
 	 };
		\addlegendentry{fail $> 10\%$}
		\addplot[blue,only marks,mark=*,mark options={solid}] coordinates {
 (27,0.098) (29,0.035)
		};
		\addlegendentry{$1\% \leq$ fail $\leq 10\%$}
		\addplot[darkgreen,only marks,mark=*,mark options={solid}] coordinates {
 (31,0.007) (33,0.005) (35,0.003) (37,0.000)
		};
		\addlegendentry{fail $< 1\%$}
		\addplot[black,dashed,mark=o,mark options={solid}] coordinates {
 (37,0.622) (39,0.253) (41,0.103) (43,0.042) (45,0.017)
		};
		\addlegendentry{Cor.~\!\ref{cor:prob_bound_reco_trig_poly}}
		\addplot[no marks,black,dotted,thick] coordinates {
 (9,1.000) (11,0.046) (13,0.018) (15,0.005) (17,0.004) (19,0.000) (21,0.000) (23,0.000)
		};
		\addlegendentry{median+R1L}
		\addplot[only marks,forget plot,red,mark=*,mark options={solid}] coordinates {
 (9,1.000)
        };
		\addplot[only marks,forget plot,blue,mark=*,mark options={solid}] coordinates {
 (11,0.046) (13,0.018)
		};
		\addplot[only marks,forget plot,darkgreen,mark=*,mark options={solid}] coordinates {
 (15,0.005) (17,0.004) (19,0.000) (21,0.000) (23,0.000)
		};
		\end{axis}
		\end{tikzpicture}
		\caption{Success rate and failure rate with respect to the number of used rank\mbox{-}1 lattices for multiple random rank-1 lattices, where $\Gamma\subset\Z^8$ is a symmetric hyperbolic cross and $I\subset \Gamma$ a weighted hyperbolic cross.
        Solid lines correspond to Algorithm~\ref{alg:compute_p_from_Gamma}, the dotted lines corresponds to Algorithm~\ref{alg:compute_p_from_Gamma} with postprocessing as discussed in Remark~\ref{rem:use_r1l_info_or_ls}, and the dashed line shows the theoretical bounds from Corollary~\ref{cor:prob_bound_reco_trig_poly}.
        }\label{fig:success_rate_L_hc_weightedhc}
\end{figure}
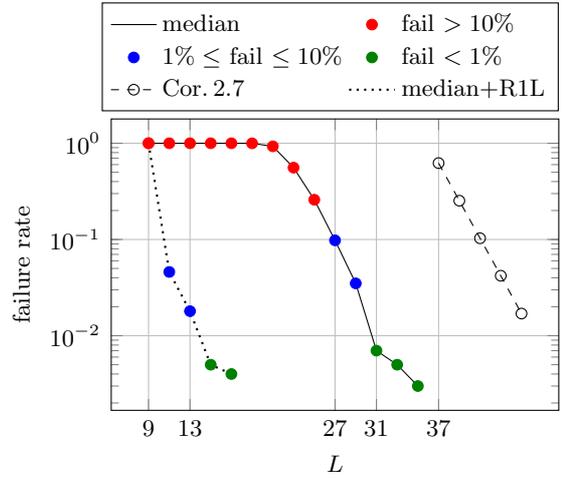

\begin{figure}
		\begin{tikzpicture}[baseline]
		\begin{axis}[
			font=\footnotesize,
			enlarge x limits=true,
			enlarge y limits=true,
			height=0.35\textwidth,
			grid=major,
			width=0.48\textwidth,
			xtick={15,19,37},
			ymin=0,ymax=4,
			xlabel={$L$},
			ylabel={max. \# (in 1\,000 runs)},
			legend style={legend cell align=left}, legend pos=south east,
			legend columns = 1,
			title={Maximal number of PFNs},
		]
		\addplot[red,mark=square*,mark options={solid}] coordinates {
 (15,4) (17,4)
		};
		\addplot[darkgreen,mark=square*,mark options={solid}] coordinates {
 (19,0) (21,0) (23,0) (25,0) (27,0) (29,0) (31,0) (33,0) (35,0) (37,0)
		};
		\end{axis}
		\end{tikzpicture}
\hfill
		\begin{tikzpicture}[baseline]
		\begin{axis}[
			font=\footnotesize,
			enlarge x limits=true,
			enlarge y limits=true,
			height=0.35\textwidth,
			grid=major,
			width=0.48\textwidth,
			xtick={19,23,25,37},
			ytick={0,10,100},
			xmin=15,
			ymin=0,ymax=100,
			xlabel={$L$},
			ylabel={max. \# (in 1\,000 runs)},
			legend style={legend cell align=left}, legend pos=south east,
			legend columns = 1,
			title={Maximal number of PFPs},
		]
		\addplot[red,no marks,mark options={solid}] coordinates {
 (15,330) (17,136) (19,58)
		};
		\addplot[black,mark=square*,mark options={solid}] coordinates {
 (19,58) (21,26) (23,14) (25,10)
		};
		\addplot[blue,mark=square*,mark options={solid}] coordinates {
 (25,10) (27,6) (29,4) (31,2) (33,2) (35,4) (37,0)
		};
		\end{axis}
		\end{tikzpicture}

		\begin{tikzpicture}[baseline]
		\begin{axis}[
			font=\footnotesize,
			enlarge x limits=true,
			enlarge y limits=true,
			height=0.35\textwidth,
			grid=major,
			width=0.48\textwidth,
			xtick={15,19,23,37},
			ytick={0,0.1,0.5,0.8,0.9,1},
			ymin=0,ymax=1,
			xlabel={$L$},
			ylabel={success rate},
			legend style={legend cell align=left}, legend pos=south east,
			legend columns = 1,
			title={Modified success rate allowing for
				$\leq$ $10$ PFPs},
            title style={at={(0.4,1)}},
		]
		\addplot[forget plot,no marks,black] coordinates {
 (15,0.000) (17,0.000) (19,0.117) (21,0.854) (23,0.995) (25,1.000) (27,1.000) (29,1.000) (31,1.000) (33,1.000) (35,1.000) (37,1.000)
		};
		\addplot[only marks,red,mark=*,mark options={solid}] coordinates {
 (15,0.000) (17,0.000) (19,0.117) (21,0.854)
		};
		\addlegendentry{success $< 90\%$}
		\addplot[darkgreen,only marks,mark=*,mark options={solid}] coordinates {
 (23,0.995) (25,1.000) (27,1.000) (29,1.000) (31,1.000) (33,1.000) (35,1.000) (37,1.000)
		};
		\addlegendentry{success $> 99\%$}
		\end{axis}
		\end{tikzpicture}
\hfill
		\begin{tikzpicture}[baseline]
		\begin{axis}[
			font=\footnotesize,
			enlarge x limits=true,
			enlarge y limits=true,
			height=0.35\textwidth,
			grid=major,
			width=0.48\textwidth,
			xtick={15,17,19,37},
			ytick={0,0.1,0.5,0.8,0.9,1},
			ymin=0,ymax=1,
			xlabel={$L$},
			ylabel={success rate},
			legend style={legend cell align=left}, legend pos=south east,
			legend columns = 1,
			title={Modified success rate allowing for
				$\leq$ $100$ PFPs},
            title style={at={(0.4,1)}},
		]
		\addplot[forget plot,no marks,black] coordinates {
 (15,0.001) (17,0.977) (19,1.000) (21,1.000) (23,1.000) (25,1.000) (27,1.000) (29,1.000) (31,1.000) (33,1.000) (35,1.000) (37,1.000)
		};
		\addplot[only marks,red,mark=*,mark options={solid}] coordinates {
 (15,0.001)
		};
		\addlegendentry{success $< 90\%$}
		\addplot[only marks,blue,mark=*,mark options={solid}] coordinates {
 (17,0.977)
		};
		\addlegendentry{$90\% \leq$ success $\leq 99\%$}
		\addplot[darkgreen,only marks,mark=*,mark options={solid}] coordinates {
 (19,1.000) (21,1.000) (23,1.000) (25,1.000) (27,1.000) (29,1.000) (31,1.000) (33,1.000) (35,1.000) (37,1.000)
		};
		\addlegendentry{success $> 99\%$}
		\end{axis}
		\end{tikzpicture}
\caption{Additional performance measures for $\Gamma$ and $I$ unweighted and weighted hyperbolic crosses, respectively, from Figure~\ref{fig:success_rate_L_hc_weightedhc} for Algorithm~\ref{alg:compute_p_from_Gamma}.\\
	{\em First row:}  Maximal numbers of aliasing frequencies within $I$ (PFNs) and aliasing frequencies in $\Gamma\setminus I$ with frequencies within $I$ (PFPs).\\
	{\em Second row:} Modified success rates with no potential false negatives, but some potential false positives allowed. }
\label{fig:success_rate_L_hc_weightedhc_detail}
\end{figure}

\subsection{Dimension-incremental sparse FFT}

We continue by numerically exploring the dimension-incremental sparse FFT method of \cite{PoVo14,KaPoVo17}, and Algorithm~\ref{algo:sfft_general} with our multiple random rank-1 lattice sampling method described in Algorithm~\ref{alg:compute_p_from_Gamma_s} in the role of Algorithm~$\operatorname{A}$.

\subsubsection{Random sparse trigonometric polynomial}
\label{sec:numerics:sfft_rand_poly}

As in \cite[Section 3.1]{PoVo14} and \cite[Section 3.1]{KaPoVo17}, we construct random multivariate trigonometric polynomials $p$ of the form~\eqref{eq:trig_pol}
with frequencies supported in the cube $\hat{G}_{N}^d=[-N,N]^d\cap\Z^d$.
For this, we choose $\vert I\vert$ frequencies $\boldk\in\hat{G}_{N}^d$ uniformly at random
and draw the corresponding Fourier coefficients $\hat{p}_\boldk\in[-1,1)+[-1,1)\mathrm{i}$, $\vert\hat{p}_\boldk\vert\geq 10^{-6}$, uniformly at random for all $\boldk\in I=\mathrm{supp}\,\hat{p}$.
For the reconstruction of the trigonometric polynomials $p$, we only assume that the search domain $\Gamma:=\hat{G}_{N}^d\supset I$.

We set the expansion parameter $N:=32$, which corresponds to $N_\Gamma=64$.
Now, we compare the results of single reconstructing rank-1 lattice sampling from \cite[Algorithm~5]{KaPoVo17}, multiple reconstructing rank-1 lattice sampling from \cite[Algorithm~4]{KaPoVo17}, and Algorithm~\ref{algo:sfft_general} combined with our multiple random rank-1 lattice approach, Algorithm~\ref{alg:compute_p_from_Gamma_s}.
Note that \cite[Algorithm~5]{KaPoVo17} also follows the framework of Algorithm~\ref{algo:sfft_general},  with a single reconstructing rank-1 lattice sampling method used in Step~3 and also taking the role of  Algorithm~$\operatorname{A}$.
Likewise, \cite[Algorithm~4]{KaPoVo17} corresponds to Algorithm~\ref{algo:sfft_general} with multiple reconstructing rank-1 lattice sampling taking the role of Algorithm~$\operatorname{A}$.
For all of these approaches, we choose the absolute threshold parameter $\theta:=10^{-12}$, the number of detection iterations $r:=1$, and the sparsity parameter $s:=\vert I\vert$.

For the combination of Algorithm~\ref{algo:sfft_general} and Algorithm~\ref{alg:compute_p_from_Gamma_s}, we work with local sparsity parameter of $s_\mathrm{local}:=2\,\vert I\vert$ and parameter $\delta:=0.9$. Motivated by the numerical results in Section~\ref{ssec:num1} and since we can tolerate up to~$|I|$ false positives for each invocation of Algorithm~\ref{alg:compute_p_from_Gamma_s},
we distinctly reduce the number $L$ of random rank-1 lattices used in Algorithm~\ref{alg:compute_p_from_Gamma_s} as compared to the choice postulated by Corollary~\ref{cor:prob_bound_reco_trig_poly}.
We only use
\begin{equation}
L:=\min\left\{n\in 2\N+1\colon n\ge \frac{1}{4}\,\frac{4c}{(c-2)\,\ln(c-1)} \, (\ln |J_t|-\,\ln\delta)\right\} \text{\;with } c:=10.33, \label{equ:L_one_half_in_numerics}
\end{equation}
which is a reduction of about 3/4 as compared to Corollary~\ref{cor:prob_bound_reco_trig_poly}.
Obviously, this choice of the parameter $L$ is not covered by our theory. However, the numerical tests below
are not affected, which impressively corroborates that the worst-case scenarios we used for the theoretical estimates hardly ever occur simultaneously.
At this point, we stress on the fact that \cite[Algorithm~4]{KaPoVo17} and \cite[Algorithm~5]{KaPoVo17} use numerically determined spatial discretizations for the
candidate sets $J_t$ in Step~2, where the approaches for constructing these spatial discretizations are also optimized with respect to the number of used sampling values.
Thus, both algorithms also already exploit the potential gaps between theory and practice.

\begin{table}[htb]
\centering
\begin{tabular}{|r|r||S[table-format=10.0]|S[table-format=9.0]|S[table-format=10.0]|}
\hline
& & \multicolumn{1}{c|}{\cite[Algorithm~5]{KaPoVo17}} & \multicolumn{1}{c|}{\cite[Algorithm~4]{KaPoVo17}} & \multicolumn{1}{c|}{Algorithm~\ref{algo:sfft_general}}  \\
& & \multicolumn{1}{c|}{(single reco.\ R1L)}& \multicolumn{1}{c|}{(multiple reco.\ R1L)} & \multicolumn{1}{c|}{using Algorithm~\ref{alg:compute_p_from_Gamma_s}}  \\
& & & & \multicolumn{1}{c|}{(random R1L)} \\
\hline
$d$\rule[0.1em]{0em}{1em} & $\vert I\vert$
& \multicolumn{1}{r|}{max.\ \#samples} & \multicolumn{1}{r|}{max.\ \#samples} & \multicolumn{1}{r|}{max.\ \#samples} \\
\hline
\hline
5\rule[0.55em]{0em}{0.55em} & 1\,000 & 6260605  &    581881  &    289914  \\
10 & 1\,000 & 20848685 &    1589349 &     649756  \\
15 & 1\,000 & 35937525 &    2599029 &     1011666  \\
20 & 1\,000 & 52361205 &    3609753 &     1373810  \\
25 & 1\,000 & 67164695 &    4621205 &     1735486  \\
30 & 1\,000 & 80660385 &    5644059 &     2097396  \\
\hline
5\rule[0.55em]{0em}{0.55em} & 10\,000 & 190618285  &    4648335  &   3321330  \\
10 & 10\,000 & 1081274675  &   15186447  &   7990386  \\
15 & 10\,000 & 1969412575  &  25662189  &   12639840  \\
20 & 10\,000 & 2935663575  & 36161887  &   17308866  \\
25 & 10\,000 & 3837073825  & 46681103  &   21958610  \\
30 & 10\,000 & 4771398905  & 57203659  &   26567030  \\
\hline
5\rule[0.55em]{0em}{0.55em} & 100\,000 & \multicolumn{1}{c|}{--}  & 33428113  &  34007204  \\ %
10 & 100\,000 & \multicolumn{1}{c|}{--}  & 143681689  & 80494280  \\ %
15 & 100\,000 & \multicolumn{1}{c|}{--}  & 250232085  & 126776914    \\ %
20 & 100\,000 & \multicolumn{1}{c|}{--} & 356857499   & 173262642    \\ %
25 & 100\,000 & \multicolumn{1}{c|}{--}  & 463174925 &  219749054   \\ %
30 & 100\,000 & \multicolumn{1}{c|}{--}  & 569711277 &  266435166  \\ %
\hline
\end{tabular}
\caption{Results for random sparse trigonometric polynomials applying \cite[Algorithm~5]{KaPoVo17}, \cite[Algorithm~4]{KaPoVo17}, and Algorithm~\ref{algo:sfft_general} using Algorithm~\ref{alg:compute_p_from_Gamma_s},
when considering frequencies within the search domain $\Gamma=\hat{G}_{32}^d$. 
The detection was successful in all considered cases and the relative $\ell_2$-errors of the Fourier coefficients near machine precision (below $2\cdot 10^{-15}$).}
\label{table:sfft:dimincr:per:numerics:rand_poly:alg2_mlfft4_mrlfft:N32}
\end{table}

\begin{sloppypar}
For sparsity $\vert I\vert\in\{1\,000, 10\,000, 100\,000\}$, we run tests for spatial dimension $d\in\{5,10,15,20,25,30\}$ applying \cite[Algorithm~5]{KaPoVo17}, \cite[Algorithm~4]{KaPoVo17}, and Algorithm~\ref{algo:sfft_general} using Algorithm~\ref{alg:compute_p_from_Gamma_s}. We only omit the case of sparsity $\vert I\vert=100\,000$ for \cite[Algorithm~5]{KaPoVo17} since this would have required quite a large number of samples and very long runtimes.
All tests are repeated $10$ times with newly chosen frequencies $\boldk\in \Gamma$ and Fourier coefficients $\hat{p}_\boldk\in \mathbb C$.
Then, for the 10 repetitions, we determine the maximum of the total number of samples.
The numerical results are displayed in Table~\ref{table:sfft:dimincr:per:numerics:rand_poly:alg2_mlfft4_mrlfft:N32}. We also determine the relative $\ell_2$-errors of the Fourier coefficients and observe that all errors are near machine precision (below $2\cdot 10^{-15}$). In particular, all frequencies in all test runs are successfully recovered.
\end{sloppypar}

For \cite[Algorithm~5]{KaPoVo17} using single reconstructing rank-1 lattices, the results are shown in the third column of Table~\ref{table:sfft:dimincr:per:numerics:rand_poly:alg2_mlfft4_mrlfft:N32}.
Moreover, the results for \cite[Algorithm~4]{KaPoVo17} using multiple reconstructing rank-1 lattices are presented in the fourth column.
We observe that we require significantly fewer samples than if we had used a $d$-dimensional FFT on a full grid, which would require $|\hat{G}_{32}^5| = 65^5 = 1\,160\,290\,625$ already in the 5-dimensional case.
Additionally, for sparsity $\vert I\vert:=10\,000$, \cite[Algorithm~4]{KaPoVo17} using multiple reconstructing rank-1 lattices required only a fraction of approximately between 1/83 and 1/11 of the samples as compared to \cite[Algorithm~5]{KaPoVo17} using single rank-1 lattices.

The results for Algorithm~\ref{algo:sfft_general} combined with the multiple random rank-1 lattice approach of Algorithm~\ref{alg:compute_p_from_Gamma_s} are presented in the fifth column of Table~\ref{table:sfft:dimincr:per:numerics:rand_poly:alg2_mlfft4_mrlfft:N32}. Here, we require only a fraction between 1/179 and 1/21 of the samples as compared to \cite[Algorithm~5]{KaPoVo17}. Moreover, we achieve at least a similar number of samples as  compared to \cite[Algorithm~4]{KaPoVo17} and are able to reduce the samples by up to 2/3.

\begin{table}[htb]
\centering
\begin{tabular}{|r|r||S[table-format=10.0]|S[table-format=9.0]|S[table-format=10.0]|}
\hline
& & \multicolumn{1}{c|}{\cite[Algorithm~5]{KaPoVo17}} & \multicolumn{1}{c|}{\cite[Algorithm~4]{KaPoVo17}} & \multicolumn{1}{c|}{Algorithm~\ref{algo:sfft_general}}  \\
& & \multicolumn{1}{c|}{(single reco.\ R1L)}& \multicolumn{1}{c|}{(multiple reco.\ R1L)} & \multicolumn{1}{c|}{using Algorithm~\ref{alg:compute_p_from_Gamma_s}}  \\
& & & & \multicolumn{1}{c|}{(random R1L)}  \\
\hline
$d$\rule[0.1em]{0em}{1em} & $\vert I\vert$
& \multicolumn{1}{r|}{max.\ \#samples} & \multicolumn{1}{r|}{max.\ \#samples} & \multicolumn{1}{r|}{max.\ \#samples} \\
\hline
\hline
   5\rule[0.55em]{0em}{0.55em} &  1000 &  81021568 &    4556111 &  372790   \\
  10 &  1000 &  218672495 &    11545999 &     842668   \\
  15 &  1000 &  345163622 &    18382277 &    1309654    \\
  20 &  1000 &  482567506 &    25400061 &    1775348    \\
  25 &  1000 &  614584563 &    32241219 &  2240656    \\
  30 &  1000 &  785132009 &    39265929 &   2712170    \\
\hline
   5\rule[0.55em]{0em}{0.55em} & 10000 & 3192096843 &   49697515 &   3745910   \\
  10 & 10000 & 11402137557 &   132626899  &   8419596    \\
  15 & 10000 & 19075479723 &   215457509  &   13069812    \\
  20 & 10000 & 27897681285 &   298481637 &    17741814    \\
  25 & 10000 & 35745041259 &   381485827 &    22373004    \\
  30 & 10000 & 43029740781 &   464049693 &  27023214    \\
\hline
 5\rule[0.55em]{0em}{0.55em} & 100\,000 & \multicolumn{1}{c|}{--} &  454666789  & 44144134  \\
10 & 100\,000 & \multicolumn{1}{c|}{--}  & 1294938567  &   100961292  \\
15 & 100\,000 & \multicolumn{1}{c|}{--}  & 2134772195  &   157779374  \\
20 & 100\,000 & \multicolumn{1}{c|}{--}  & 2974838029  &   214596532  \\
25 & 100\,000 & \multicolumn{1}{c|}{--}  & 3815598519  &   271616236  \\
30 & 100\,000 & \multicolumn{1}{c|}{--}  & 4655793605  &   328233342  \\
\hline
\end{tabular}
\caption{Results for random sparse trigonometric polynomials applying \cite[Algorithm~5]{KaPoVo17}, \cite[Algorithm~4]{KaPoVo17}, and Algorithm~\ref{algo:sfft_general} using Algorithm~\ref{alg:compute_p_from_Gamma_s},
when considering frequencies within the search domain $\Gamma=\hat{G}_{256}^d$.
The detection was successful in all considered cases and the relative $\ell_2$-errors of the Fourier coefficients near machine precision (below $2\cdot 10^{-15}$).}
\label{table:sfft:dimincr:per:numerics:rand_poly:alg2_mlfft4_mrlfft:N256}
\end{table}

Furthermore, we also re-ran all the test for an increased expansion parameter $N=256$, which corresponds to $N_\Gamma=512$. The numerical results we obtained are shown in Table~\ref{table:sfft:dimincr:per:numerics:rand_poly:alg2_mlfft4_mrlfft:N256}; the columns have the same meaning as in Table~\ref{table:sfft:dimincr:per:numerics:rand_poly:alg2_mlfft4_mrlfft:N32}.
For this scenario, we observe that Algorithm~\ref{algo:sfft_general} using the multiple random \hbox{rank-1} lattice approach from Algorithm~\ref{alg:compute_p_from_Gamma_s} requires only a fraction between 1/17 and 1/10 of the number of samples of \cite[Algorithm~4]{KaPoVo17} using multiple reconstructing rank-1 lattices as well as only a fraction between 1/1590 and 1/217 of the number of samples of \cite[Algorithm~5]{KaPoVo17}. The relative $\ell_2$-errors are still near machine precision for all the  settings and methods considered (below $2\cdot 10^{-15}$).

Comparing the results in Table~\ref{table:sfft:dimincr:per:numerics:rand_poly:alg2_mlfft4_mrlfft:N256} for expansion parameter $N=256$ with the results in Table~\ref{table:sfft:dimincr:per:numerics:rand_poly:alg2_mlfft4_mrlfft:N32} for $N=32$,
we observe a significantly higher reduction in the number of samples for the combination of Algorithm~\ref{algo:sfft_general} and Algorithm~\ref{alg:compute_p_from_Gamma_s} in the case $N=256$. This observation confirms the theoretical results that we are able to reduce the factor $N_\Gamma$ in the sample complexity of \cite[Algorithm~5]{KaPoVo17} and \cite[Algorithm~4]{KaPoVo17} to $\log N_\Gamma$ by employing Algorithm~\ref{alg:compute_p_from_Gamma_s}, see also Table~\ref{tab:complexities_comparison}.

\subsubsection{Approximation of tensor-product function by trigonometric polynomials}
\label{sec:numerics:sfft_fct}
In the following, we will demonstrate that our method also works for functions that are only approximately sparse.
 For that we consider the multivariate periodic test function $f\colon\T^{10}\rightarrow\R$,
\begin{equation} \label{equ:f:10}
f(\boldx):=\prod_{t\in\{1,3,8\}}N_2(x_t) + \prod_{t\in\{2,5,6,10\}}N_4(x_t) + \prod_{t\in\{4,7,9\}}N_6(x_t),
\end{equation}
from \cite[Section~3.3]{PoVo14} and \cite[Section~3.3]{KaPoVo17} 
where $N_m:\T\rightarrow\R$ is the B-Spline of order $m\in\N$,
$$N_m(x) := C_m \sum_{k\in\Z} \operatorname{sinc}\left(\frac{\pi}{m}k\right)^m (-1)^k \,\mathrm{e}^{2\pi\mathrm{i}kx},$$
with a constant $C_m>0$ such that $\Vert N_m \vert L_2(\T)\Vert=1$.

Each of the three summands has infinitely many non-zero Fourier coefficients, but they are exhibiting a decay pattern described by a hyperbolic cross (different for each of the summands). Thus we expect that the  function is well approximated by trigonometric polynomials, namely the one corresponding to the Fourier coefficients indexed by a union of three hyperbolic crosses, each corresponding to significant coefficients of one of the summands.

We aim to demonstrate that our method allows to efficiently find such an approximation of the function $f$ by a multivariate trigonometric polynomial $p$.
More precisely, we apply the dimension-incremental approaches already discussed in the previous examples to  determine a frequency index set
$I=I^{(1,\ldots,10)}\subset\Gamma:=\hat{G}_N^{10}$ and to compute
approximated Fourier coefficients $\tilde{\hat{p}}_\boldk$, $\boldk\in I$. %
As explained, an adaquate choice for the resulting frequency index sets $I$ is given by the union of three sets of frequencies corresponding to the significant coefficients of the three summands. In our example, these are
a three-dimensional symmetric hyperbolic cross in the dimensions $1,3,8$,
a four-dimensional symmetric hyperbolic cross in the dimensions $2,5,6,10$,
and a three-dimensional symmetric hyperbolic cross in the dimensions $4,7,9$.

All tests are performed 10 times and the relative $L_2(\T^{10})$ approximation error
$$
 \Vert f-\tilde{S}_I f\vert L_2(\T^{10})\Vert / \Vert f\vert L_2(\T^{10})\Vert
 =
 \sqrt{\Vert f\vert L_2(\T^{10})\Vert^2 - \sum_{\boldk\in I}\vert\hat{f}_\boldk\vert^2 + \sum_{\boldk\in I}\vert\tilde{\hat{p}}_\boldk-\hat{f}_\boldk\vert^2} / \Vert f\vert L_2(\T^{10})\Vert
$$
is computed,
where $p=\tilde{S}_I f:=\sum_{\boldk\in I} \tilde{\hat{p}}_\boldk \,\mathrm{e}^{2\pi\mathrm{i}\boldk\cdot\circ}$.

We set the expansion parameter $N:=16,32,64$
and we use the full grids $\Gamma:=\hat{G}_N^{10}$ as search space.
Moreover, we set the number of detection iterations $r:=5$.
The used sparsity input parameters $s$ and $s_\mathrm{local}=2s$ are specified in column 2 of Table~\ref{table:numerics:fct:s_sparse:a2r1l_mlfft6_mrlfft}. Furthermore, the results
of \cite[Algorithm~5]{KaPoVo17} based on single reconstructing rank-1 lattices, of \cite[Algorithm~4]{KaPoVo17} based on multiple reconstructing rank-1 lattices, and of Algorithm~\ref{algo:sfft_general} combined with our multiple random rank-1 lattice approach in Algorithm~\ref{alg:compute_p_from_Gamma_s}
are shown in columns 3--4, 5--6, and 7--8 of Table~\ref{table:numerics:fct:s_sparse:a2r1l_mlfft6_mrlfft}, respectively.
For the combination of Algorithm~\ref{algo:sfft_general} and~\ref{alg:compute_p_from_Gamma_s}, we set the parameter $\delta:=0.999$ and the number $L$ of random rank-1 lattices to~\eqref{equ:L_one_half_in_numerics}, which again corresponds to a reduction of about 3/4 compared to the theoretical predictions of Corollary~\ref{cor:prob_bound_reco_trig_poly}.

The column ``max.\ rel.\ $L_2$-error'' contains the maximum of
the relative $L_2(\T^{10})$ approximation errors $\Vert f-\tilde{S}_I f\vert L_2(\T^{10})\Vert / \Vert f\vert L_2(\T^{10})\Vert$ of the 10 test runs.
The remaining columns have the same meaning as described in Section~\ref{sec:numerics:sfft_rand_poly}.
We observe that for increasing sparsity parameter, the number of samples increases
while the relative $L_2(\T^{10})$ approximation error decreases.

Moreoever, we observe that \cite[Algorithm~5]{KaPoVo17} and Algorithm~\ref{algo:sfft_general} combined with Algorithm~\ref{alg:compute_p_from_Gamma_s} yield similar relative $L_2(\T^{10})$ approximation errors, whereas \cite[Algorithm~4]{KaPoVo17} produces slightly higher errors. Furthermore, Algorithm~\ref{algo:sfft_general} combined with Algorithm~\ref{alg:compute_p_from_Gamma_s} requires the fewest samples, a fraction of between 1/9 and 1/2 of the number of samples required by \cite[Algorithm~4]{KaPoVo17} and a fraction between 1/148 and 1/17 of the number of samples required by \cite[Algorithm~5]{KaPoVo17}.
\newline
For instance in the case $N=64$ and $s=5\,000$, the maximal total number of samples for \cite[Algorithm~5]{KaPoVo17} (computed over 10 test runs) was about 2.3 billion, for \cite[Algorithm~4]{KaPoVo17} about 159 million samples, and for Algorithm~\ref{algo:sfft_general} using Algorithm~\ref{alg:compute_p_from_Gamma_s} about 19 million samples, while achieving comparable errors.

\begin{table}[htb]
\centering
\begin{small}
\begin{tabular}{|r|r||S[table-format=10.0]|r||S[table-format=9.0]|r||S[table-format=9.0]|r|}
\hline
& & \multicolumn{2}{c||}{\cite[Algorithm~5]{KaPoVo17}} & \multicolumn{2}{c||}{\cite[Algorithm~4]{KaPoVo17}} & \multicolumn{2}{c|}{Algorithm~\ref{algo:sfft_general}} \\
& & \multicolumn{2}{c||}{(single reco.\ R1L)} & \multicolumn{2}{c||}{(multiple reco.\ R1L)} & \multicolumn{2}{c|}{using Algorithm~\ref{alg:compute_p_from_Gamma_s}}  \\
& & \multicolumn{2}{c||}{} & \multicolumn{2}{c||}{} & \multicolumn{2}{c|}{(random R1L)}  \\
    & & \multicolumn{1}{r|}{max.} & max.\ rel. & \multicolumn{1}{r|}{max.} & max.\ rel. & \multicolumn{1}{r|}{max.} & max.\ rel. \\
$N$ & $s$ & \multicolumn{1}{r|}{\#samples} & $L_2$-error & \multicolumn{1}{r|}{\#samples} & $L_2$-error & \multicolumn{1}{r|}{\#samples} & $L_2$-error \\
\hline
16\rule[0.55em]{0em}{0.55em} &  1\,000 & 50405091 & 1.2e-02 & 8094293 & 1.3e-02 &   2903576 & 1.2e-02 \\
16 &  2\,000 &  128362707 & 4.3e-03 &  15449367 & 5.0e-03 &   5813898 & 4.1e-03 \\
16 &  3\,000 &  222662847 & 3.4e-03 &  22285615 & 4.1e-03 &   9643162 & 3.1e-03 \\
\hline
32\rule[0.55em]{0em}{0.55em} &  1\,000 &  84241365 & 1.2e-02 & 16623363 & 1.3e-02 &   2905176 & 1.2e-02 \\
32 &  2\,000 &  265019105 & 3.4e-03 & 32226797 & 8.8e-03 &   6683344 & 3.4e-03 \\
32 &  3\,000 &  565847035 & 1.7e-03 & 46948295 & 2.2e-03 &   10637178 & 1.7e-03 \\
32 &  4\,000 &  767068055 & 1.4e-03 & 59623261 & 1.6e-03 &   14175646 & 1.3e-03 \\
\hline
64\rule[0.55em]{0em}{0.55em} &  1\,000 &      242940304 & 1.2e-02 &  34236011 & 4.5e-01 &   3540792 & 1.2e-02 \\
64 &  2\,000 &      563767023 & 3.4e-03 &  66914731 & 3.5e-03 &   7504972 & 3.4e-03 \\
64 &  3\,000 &   946215129 & 1.6e-03 &  98926095 & 1.8e-03 &   11272744 & 1.6e-03 \\
64 &  4\,000 &   1266499683 & 9.8e-04 &  130967911 & 1.1e-03 &   15005506 & 9.8e-04 \\
64 &  5\,000 &    2285094003 & 7.1e-04 &  158693803 & 9.6e-04 &   18772634 & 7.0e-04 \\
64 & \!\!\! 10\,000 &    5577419619 & 4.3e-04 &  298886055 & 5.1e-04 & 37534358 & 3.9e-04 \\
\hline
\end{tabular}
\end{small}
\caption{Results for function $f\colon\T^{10}\rightarrow\R$ from \eqref{equ:f:10} when limiting the number of detected frequencies, $s_\mathrm{local}:=2s$.}
\label{table:numerics:fct:s_sparse:a2r1l_mlfft6_mrlfft}
\end{table}

\subsection*{Acknowledgements}
L. K\"ammerer gratefully acknowledges funding
by the Deutsche Forschungsgemeinschaft (DFG, German Research Foundation, project
number 380648269). F.\ Krahmer gratefully acknowledges funding
by the Deutsche Forschungsgemeinschaft (DFG, German Research Foundation, Emmy Noether Junior Research Group 4512/1-1).
T.\ Volkmer gratefully acknowledges partial funding by the S\"achsische Aufbaubank -- F\"orderbank -- (SAB) 100378180.

\begin{small}

\end{small}

\end{document}